
\documentstyle{amsppt}
\magnification=1200
\hoffset=-0.5pc
\vsize=57.2truepc
\hsize=38.5truepc
\nologo
\spaceskip=.5em plus.25em minus.20em

\define\de{\partial}
\define\armcusgo{1}
\define\atibottw{2}
\define\goldmone{3}
\define\goldmtwo{4}
 \define\poisson{5}
 \define\souriau{6}
 \define\singula{7}
\define\singulat{8}
\define\topology{9}
  \define\smooth{10}
\define\singulth{11}
 \define\locpois{12}
 \define\modusym{13}
\define\huebjeff{14}
\define\karshone{15}
\define\marswein{16}
\define\narashed{17}
\define\naramntw{18}
\define\naramnth{19}
\define\seshaone{20}
\define\sjamlerm{21}
\define\weinstwe{22}
\define\whitnone{23}
\noindent
dg-ga/9411009
\smallskip\noindent
Annales de l'Institut Fourier, to appear.
\bigskip
\topmatter
\title
Poisson structures on  certain moduli spaces\\
for bundles on a surface
\endtitle
\author Johannes Huebschmann{\dag}
\endauthor
\affil
Max Planck Institut f\"ur Mathematik
\\
Gottfried Claren-Str. 26
\\
D-53 225 BONN
\\
huebschm\@mpim-bonn.mpg.de
\endaffil

\keywords{Geometry of principal bundles,
singularities of smooth mappings,
symplectic reduction with singularities,
Yang-Mills connections,
stratified symplectic space,
Poisson structure,
geometry of moduli spaces,
representation spaces,
categorical quotient,
geometric invariant theory,
moduli of vector bundles}
\endkeywords
\subjclass{14D20, 32G13, 32S60, 58C27, 58D27, 58E15,  81T13}
\endsubjclass

\thanks{{\dag} The author carried out this work in the framework
of the VBAC research group of Europroj.}
\endthanks

\abstract{
Let $\Sigma$ be a closed  surface, $G$ a compact Lie group, with Lie algebra
$g$, and $\xi \colon P \to \Sigma$ a principal $G$-bundle. In earlier work we
have shown that the moduli space $N(\xi)$ of central Yang-Mills connections,
with reference to appropriate additional data, is stratified by smooth
symplectic manifolds and  that  the holonomy yields a homeomorphism  from
$N(\xi)$ onto a certain representation space $\roman{Rep}_{\xi}(\Gamma,G)$,
in fact a diffeomorphism, with reference to suitable smooth structures
$C^{\infty}(N(\xi))$ and $C^{\infty}\left(\roman{Rep}_{\xi}(\Gamma,G)\right)$,
where $\Gamma$ denotes  the universal central extension of the fundamental
group of $\Sigma$. Given a coadjoint action invariant symmetric bilinear form
on $g^*$, we construct here Poisson structures on $C^{\infty}(N(\xi))$ and
$C^{\infty}\left(\roman{Rep}_{\xi}(\Gamma,G)\right)$ in such a way that the
mentioned diffeomorphism identifies them. When the form on $g^*$ is
non-degenerate the Poisson structures are  compatible with the
stratifications where $\roman{Rep}_{\xi}(\Gamma,G)$ is endowed with the
corresponding stratification and, furthermore,  yield structures of a
{\it stratified  symplectic space\/}, preserved by the induced action of the
mapping class  group of $\Sigma$.}
\endabstract
\endtopmatter
\document
\rightheadtext{Poisson structures on moduli spaces}
\leftheadtext{Johannes Huebschmann}

\beginsection Introduction

Let $X$ be a stratified space.
A {\it stratified symplectic structure\/}
on $X$
in the sense of {\smc Sjamaar-Lerman\/}~\cite\sjamlerm\
is a  Poisson algebra
of continuous functions
on $X$ which, on each stratum, restricts to a symplectic
Poisson algebra
of smooth functions.
In the present paper
we construct
such structures and related ones
on certain moduli spaces.
It is the fifth of a series of papers
about
a program
revealing
the structure of these moduli spaces
by means of the symplectic
or more generally Poisson
geometry
of certain related classical constrained systems.
Its predecessors are
\cite\singula~--~\cite\smooth,
and it will be followed by \cite\singulth\ and \cite\locpois.
\smallskip
We explain briefly the moduli spaces:
Let
$\Sigma$ be a closed surface,
$G$ a compact
Lie group,
not necessarily connected,
with Lie algebra $g$,
and $\xi \colon P \to \Sigma$
a principal $G$-bundle,
having a connected total space $P$.
Then a choice
of Riemannian metric on $\Sigma$ and
{\it orthogonal structure\/}
on $g$,
that is, adjoint action invariant scalar product,
determines a Yang-Mills functional
on the space $\Cal A(\xi)$ of connections
on $\xi$;
see \cite\atibottw\  to which we refer for background and notation.
We assume throughout that
solutions of the corresponding Yang-Mills equations
exist; this will be so for example
when $G$ is connected, cf. \cite\atibottw.
Then the {\it moduli space\/} $N(\xi)$
of gauge equivalence classes of
central Yang-Mills
connections
is non empty; it is a compact space,
including as special cases
moduli spaces of flat connections and
the {\smc Narasimhan-Seshadri}-moduli spaces \cite\narashed\
of semi stable holomorphic vector bundles.
In this paper we shall complete the proof of the
following.

\proclaim{Theorem 1}
The decomposition of $N(\xi)$
according to orbit types of central Yang-Mills connections
is a stratification (in the strong sense)
and the data
determine a stratified symplectic structure
$(C^{\infty}(N(\xi)),\{\cdot,\cdot\})$ for it.
\endproclaim

\smallskip
The  Poisson structure goes {\it beyond\/}
usual symplectic geometry;
it encapsulates the {\it mutual positions\/}
of the {\it symplectic structures\/}
on the strata.
\smallskip
In \cite\singulat\ we have shown that
$N(\xi)$ is stratified by smooth  manifolds
and that the strata inherit symplectic structures from the data.
In \cite\topology\
we constructed a
homeomorphism
$\rho_{\flat}$,
referred to as {\it Wilson loop mapping\/},
from $N(\xi)$
onto  a certain
representation space
$\roman {Rep}_{\xi}(\Gamma,G)$
for the universal central extension $\Gamma$ of the
fundamental group $\pi$ of $\Sigma$.
While the space $N(\xi)$
{\it depends\/}
on the choices of Riemannian metric
on $\Sigma$
the space
$\roman{Rep}_{\xi}(\Gamma,G)$
does not.
In
\cite\smooth\
we constructed
smooth structures
$C^{\infty}(N(\xi))$
and $C^{\infty}\left(\roman {Rep}_{\xi}(\Gamma,G)\right)$
on these spaces,
and we have shown that
$\rho_{\flat}$
is a diffeomorphism with respect to these structures.
In the present paper,
proceeding somewhat more generally than needed
for the proof of the above theorem,
we construct Poisson structures on the algebras
$C^{\infty}(N(\xi))$
and
$C^{\infty}\left(\roman{Rep}_{\xi}(\Gamma,G)\right)$
involving as additional ingredient
a coadjoint action invariant symmetric bilinear form
on $g^*$,
not necessarily positive definite
nor non-degenerate,
in such a way that
the Wilson loop mapping
identifies the Poisson structures.
When the bilinear form
on $g^*$
is positive definite
we obtain the Poisson algebra in the above Theorem.
We now explain informally the Poisson brackets.
\smallskip
By construction,
the space $\roman{Rep}_{\xi}(\Gamma,G)$
is the quotient $\roman{Hom}_{\xi}(\Gamma,G)\big / G$
of a certain space
of homomorphisms
$\roman{Hom}_{\xi}(\Gamma,G)$
of $\Gamma$ into $G$
determined by $\xi$; see Section 1 below for details.
For
$\phi \in \roman{Hom}_{\xi}(\Gamma,G)$, we denote by
$g^*_{\phi}$
the dual $g^*$
of the Lie algebra $g$, made into a
$\pi$-module via $\phi$
and the coadjoint action.
For every $[\phi] \in \roman{Rep}_{\xi}(\Gamma,G)$,
a choice of representative
$\phi \in \roman{Hom}_{\xi}(\Gamma,G)$
induces a linear map
$\lambda^*_{\phi}$
from the real vector space
$\Omega_{[\phi]}\roman{Rep}_{\xi}(\Gamma,G)$
of differentials at $[\phi]$,
with reference to
$C^{\infty}\left(\roman{Rep}_{\xi}(\Gamma,G)\right)$,
into the
first homology group
$\roman H_1(\pi,g^*_{\phi})$
of $\pi$ with coefficients in
$g^*_{\phi}$, and
$\lambda^*_{\phi}$
is
an isomorphism if and only if
$[\phi]$ is a non-singular point of $\roman{Rep}_{\xi}(\Gamma,G)$,
see (1.16) below.
In view of \cite\smooth\ (7.10),
$\lambda^*_{\phi}$ is independent of the choice of $\phi$ in the sense that,
given $x \in G$,
the corresponding
linear map
$\lambda^*_{x\phi}$
from
$\Omega_{[\phi]}\roman{Rep}_{\xi}(\Gamma,G)$
to
$\roman H_1(\pi,g^*_{x\phi})$
equals the composite
of $\lambda^*_{\phi}$ with the isomorphism
$\roman{Ad^*}(x)_{\flat}$
induced by $x$, cf. (1.16) below.
A coadjoint action invariant
symmetric bilinear form
$
\langle\cdot,\cdot\rangle
$
on $g^*$,
not necessarily
positive definite nor
non-degenerate,
then gives rise to,
for every $\phi \in \roman{Hom}_{\xi}(\Gamma,G)$,
an intersection pairing
$\langle\cdot,\cdot\rangle_{\phi}$
on
$\roman H_1(\pi,g^*_{\phi})$,
and
the corresponding Poisson bracket $\{\cdot,\cdot\}$
on
$C^{\infty}\left(\roman{Rep}_{\xi}(\Gamma,G)\right)$
will then satisfy the formula
$$
\{f,h\}[\phi]
=
\langle \lambda^*_{\phi}(df[\phi]),\lambda^*_{\phi}(dh[\phi])\rangle_{\phi}
\tag0.1
$$
where
$f,h \in C^{\infty}\left(\roman{Rep}_{\xi}(\Gamma,G)\right)$ and where
$\phi \in \roman{Hom}_{\xi}(\Gamma,G)$
is a
representative
of the point $[\phi] \in \roman{Rep}_{\xi}(\Gamma,G)$.
As a set function, the
bracket $\{\cdot,\cdot\}$  is determined by (0.1).
This formula is {\it intrinsic\/} in the sense that it
does {\it not\/} involve choices
except that of the representative
$\phi$ which has been taken care of
already by the discussion of the dependence of
$\lambda^*_{\phi}$ on the choice of $\phi$.
It then remains to prove that
(i)
for every
$f,h \in C^{\infty}\left(\roman{Rep}_{\xi}(\Gamma,G)\right)$,
the bracket
$\{f,h\}$ is an element
of $C^{\infty}\left(\roman{Rep}_{\xi}(\Gamma,G)\right)$
and, the Leibniz rule and skew symmetry
being obviously true for
the bracket $\{\cdot,\cdot\}$, that
(ii) this bracket satisfies the Jacobi identity.
An appropriately  reworded statement will be given in (2.1) below.
Proofs will then be given in Section 2 except that Section 3
is devoted to the Jacobi identity.
Our construction is entirely finite dimensional except that
the {\it verification of the Jacobi identity\/}
involves
(i) the smooth open connected
and dense stratum whose existence has been established
in \cite\singulat\
and, furthermore, (ii) the local model constructed in \cite\singula,
for
the special case where the bilinear
form on $g^*$ is non-degenerate.
The general case of
an arbitrary symmetric bilinear
2-form on $g^*$
is
then handled by relating it to that of a certain associated
non-degenerate 2-form. See Section 3 below for details.
For a non-degenerate
2-form on $g^*$,
the
Jacobi identity on the smooth open connected and dense
stratum can also
be settled by the finite dimensional techniques in \cite\weinstwe.
\smallskip
The intrinsic description (0.1) of the Poisson
structure
has the following consequence
a proof of which will be given at the end of Section 2.

\proclaim{Theorem 2}
The induced action
of the mapping class group of $\Sigma$
respects the Poisson structure.
More precisely, its subgroup
of orientation
preserving elements
preserves the Poisson bracket
on $\roman{Rep}_{\xi}(\Gamma,G)$
whereas
the orientation
reversing elements
yield
Poisson bracket preserving diffeomorphisms
from
$\roman{Rep}_{\xi}(\Gamma,G)$
to
$\roman{Rep}_{-\xi}(\Gamma,G)$
where $-\xi$ refers to the (topologically)
\lq\lq opposite\rq\rq\  bundle
(which may coincide with $\xi$).
\endproclaim

\smallskip
When the bilinear
form on $g^*$ is non-degenerate
but not necessarily positive definite
the resulting Poisson structure is {\it symplectic\/} in the sense
that its Casimir elements
are the constants only
and
in fact then yields a structure
of a {\it stratified symplectic space\/},
cf. Sections  4 below;
this then completes the proof of
a somewhat more general result than Theorem 1;
see (4.1) below for details.
In Section 5 we indicate how the
twist flows constructed in \cite\goldmtwo\
on the top stratum can be extended to the whole space
as derivations of the smooth structure.
Section 1 below is preparatory in character.
\smallskip
It has been known for a while,
cf. e.~g. {\smc Narasimhan-Seshadri\/}~\cite\narashed,
{\smc Atiyah-Bott\/}~\cite\atibottw,
{\smc Goldman\/}~\cite\goldmone,
that an orthogonal structure on the Lie algebra
$g$ gives rise to a symplectic structure on
a certain {\it non-singular part\/}
of spaces of the kind
$N(\xi)$ and $\roman{Rep}_{\xi}(\Gamma,G)$.
However, in general, these spaces come with singularites,
and our Poisson structures include these singularities.
One of the chief results of our earlier paper,
\cite\smooth\ (6.2) and (6.3), says that,
locally, a space of the kind
$N(\xi)$ and $\roman{Rep}_{\xi}(\Gamma,G)$
looks like the reduced space
of a momentum mapping
for a representation of a compact Lie group varying over the
space. The result of the present paper
says that this local picture is available
even for a {\it globally\/} defined Poisson structure
which, in the local model, then amounts to
the {\smc Arms-Cushman-Gotay} Poisson structure
\cite\armcusgo\
on the reduced space for a representation
of a compact Lie group.
The local model, with the Poisson structure included, is made precise
in (4.3) below.
We hope to prove elsewhere
that a suitable holomorphic quantization
of such a globally defined Poisson structure
then yields
a finite dimensional
complex vector space.
\smallskip
An illustration of our result
is worked out in our paper \cite\locpois:
For $G=\roman{SU}(2)$,
in another guise, the moduli space $N(\xi)$ is that of semi stable
holomorphic vector bundles on $\Sigma$
(with reference to a choice of holomorphic structure)
of rank 2, degree 0, and trivial determinant.
This space and related ones have been studied extensively in the literature
\cite\naramntw\  -- \cite\seshaone.
In particular,
for genus $\ell \geq 2$,
the complement $\Cal K$ of the top stratum
is known to be the {\it Kummer\/} variety
of $\Sigma$ associated with its Jacobien $J$
and the canonical involution thereupon.
In \cite\locpois\ we prove the following.

\proclaim{Theorem 3}
When $\Sigma$ has genus $\ell \geq 2$,
the
Poisson algebra $(C^{\infty}(N(\xi)),\{\cdot,\cdot\})$
detects the Kummer variety
$\Cal K$
in $N(\xi)$
together with its
$2^{2 \ell}$ double points.
More precisely,
$\Cal K$
consists of the points of $N(\xi)$ where the rank of
the Poisson structure
is not maximal, the double points
being those where the rank is zero.
\endproclaim

In particular,
when $\Sigma$ has genus two,
the space $N(\xi)$  equals complex projective
3-space
and $\Cal K$ is the Kummer surface
associated with the Jacobien of $\Sigma$, cf.
{\smc Narasimhan-Ramanan}~\cite\naramntw.
In the literature,
this case has been considered somewhat special
since as a {\it space\/}
$N(\xi)$ is then actually smooth.
However, from our point of view, there is
{\it no\/}
exception.
As a {\it stratified  symplectic\/} space,
$N(\xi)$ still has singularities,
our algebra
$C^{\infty}(N(\xi))$ is {\it not\/}
that of smooth functions in the ordinary sense,
and the Kummer surface $\Cal K$
is the complement of the top stratum
and hence still precisely the singular locus
in the sense of stratified symplectic space;
in particular,
the symplectic structure
on the top stratum does {\it not\/}
extend to the whole space.
It is interesting to observe that the symplectic stratification,
that is, the one used exclusively in our approach,
is {\it finer\/}
than the standard complex analytic one on complex
projective 3-space.
\smallskip
Recent work of Jeffrey and the author
\cite\huebjeff\
shows that the moduli spaces
can be obtained by finite dimensional reduction
from certain \lq\lq extended moduli spaces\rq\rq.
The theorem of Sjamaar and Lerman
can then be applied to obtain the stratified
symplectic structures
on the moduli space.
The approach of the current paper is elementary
but less elegant than that in
\cite\huebjeff;
yet it has its advantages:
It yields at once the intrinsic formula
(0.1) above from which
the compatibility of the structure with the
action of the mapping class group is deduced.
\smallskip
I am indebted to
M. S. Narasimhan and
A. Weinstein for their interest in the present
work and for discussions at various stages of the project.

\medskip\noindent{\bf 1. Differentials}\smallskip\noindent
\define\rel{r}
Pick a base point
$Q$ of $\Sigma$ and
consider the standard presentation
$$
\Cal P  = \big\langle x_1,y_1,\dots, x_\ell,y_\ell; r\big\rangle,
\quad
r = \left[x_1,y_1\right] \cdot
\dots
\cdot
\left[x_\ell,y_\ell\right] ,
\tag1.1
$$
of the fundamental group
$\pi=\pi_1(\Sigma,Q)$,
the number $\ell$ being
the genus of  $\Sigma$;
we denote by $F$ the free group on
$x_1,y_1,\dots, x_\ell,y_\ell$ and by $N$ the normal closure
of $r$ in $F$ so that $\pi = F/N$.
For an arbitrary commutative ground ring $R$,
the presentation $\Cal P$
yields
the object
$$
\widehat{\roman R(\Cal P)}
\colon
RF
@<{\de_1^F}<<
RF\left[x_1,y_1,\dots,x_\ell,y_\ell\right]
@<{\de_2^F}<<
RF\left[r \right]
\tag1.2
$$
cf. \cite\smooth\ (5.2).
Here
$RF\left[x_1,y_1,\dots,x_\ell,y_\ell\right]$ and
$RF\left[r \right]$ refer to the free
right $RF$-modules
having $x_1,y_1,\dots,x_\ell,y_\ell$ and
$r$ as bases, respectively,
and,
the elements of these modules
being viewed as column vectors,
the operators $\de_*^F$ are given by
$$
\de_2^F =
\left[
\frac {\de r}{\de x_1},
\cdots,
\frac {\de r}{\de y_{\ell}}
\right]^{\roman{t}}\colon  RF \to (RF)^{2\ell} ,
\quad
\de_1^F
=
\left[
1-x_1,
\cdots,
1-y_\ell
\right] \colon (RF)^{2\ell} \to RF,
$$
where ${}^{\roman{t}}$ refers to the transpose of a vector.
Given a left $RF$-module $V$,
with structure map
$\chi$ from  $F$ to $\roman{Aut}(V)$,
application of the functor $-\otimes _F V$
to
$\widehat{\roman R(\Cal P)}$
yields the sequence
${\left(\widehat{\roman R(\Cal P)}\otimes_{RF} V,\partial_*^{\chi} \right)}$
which,
in view of the obvious identifications
$\widehat{\roman R_2(\Cal P)}\otimes_{RF} V = V^m$,
$\widehat{\roman R_1(\Cal P)}\otimes_{RF} V = V^n$,
and $\widehat{\roman R_0(\Cal P)}\otimes_{RF} V = V$,
looks like
$$
V
@<{\partial_1^{\chi}}<<
V^{2\ell}
@<{\partial_2^{\chi}}<<
V.
\tag1.3
$$
Here the operators
$\partial_*^{\chi}$ depend on the
$RF$-module structure on $V$
whence the notation.
\smallskip
Modulo $N$, (1.2)
yields the
free resolution
$$
\roman R(\Cal P)
\colon
R\pi
@<{\de_1}<<
R\pi\left[x_1,y_1,\dots,x_\ell,y_\ell\right]
@<{\de_2}<<
R\pi\left[r \right]
\tag1.4
$$
of the ground ring $R$, viewed as a trivial $R\pi$-module,
in the category of {\it right\/} $R\pi$-modules.
Thus when
the left $RF$-module structure
$\chi$
on $V$
factors through a left $R\pi$-module structure on $V$,
the sequence (1.3) is a chain complex
$$
\bold C(\Cal P,V)\colon
C_0(\Cal P,V)
@<{\partial_1^{\chi}}<<
C_1(\Cal P,V)
@<{\partial_2^{\chi}}<<
C_2(\Cal P,V)
\tag1.5
$$
computing the homology groups of
$\pi$ with coefficients in $V$.
\smallskip
We now take
$R=\bold R$, the reals, and
$V=g^*$, the dual of the Lie algebra $g$ of $G$,
with the corresponding structure of a
{\it left\/} $G$-module
given by the coadjoint action,
that is,
for a linear form $u$ on $g$  and  $x \in G$,
$$
(x u) = u \circ \roman{Ad}(x^{-1})\colon g \to \bold R.
\tag1.6
$$
Given
a homomorphism $\phi$ from  $F$ to $G$, we write
$g^*_{\phi}$
for $g^*$, viewed as a left
$F$-module
via $\phi$ and the coadjoint action of $G$ on $g^*$;
the sequence (1.3) then looks like
$$
g^*
@<{\partial_1^{\phi}}<<
\left(g^*\right)^{2\ell}
@<{\partial_2^{\phi}}<<
g^*;
\tag1.7
$$
when
$\phi$
is a homomorphism
from $F$ to $G$
so that each $\phi(r_j)$
lies in the centre of $G$,
the left $F$-module $g^*_{\phi}$
inherits a structure of a left $\pi$-module
which we still denote by
$g^*_{\phi}$,
and (1.7) computes the homology $\roman H_*(\pi, g^*_{\phi})$.
\smallskip
For $w \in G$, the operation of left translation
from $g$ to $\roman T_wG$
will be written
$\roman L_w$.
The assignment
to $\chi \in \roman{Hom}(F,G)$
of $(\chi(x_1),\dots,\chi(y_{\ell})) \in G^{2\ell}$
identifies $\roman{Hom}(F,G)$
with $G^{2\ell}$.
The
homomorphism
$\chi$
being viewed as the point
${
\bold w
=\left(\chi(x_1),\dots,\chi(y_{\ell})\right)}$
of $G^{2\ell}$,
its operation
of {\it left translation\/} is then
an isomorphism
$\roman L_{\chi}$ from
$g^{2\ell}$ to
$\roman T_{\chi}\roman{Hom}(F,G)$.
At $\chi\in\roman{Hom}(F,G)$,
with reference to the usual smooth structure,
the vector space of differentials
$\roman T^*_{\chi}\roman{Hom}(F,G)$
is just the usual
cotangent space,
and the dual of
$\roman L_{\chi}$ is an isomorphism
$\roman L^*_{\chi}$ from
$\roman T^*_{\chi}\roman{Hom}(F,G)$
to
$(g^*)^{2\ell}$.
Henceforth we confuse in notation
the relator $r$
with its word map
from $\roman{Hom}(F,G)$ to $G$; it is given by the assignment
to
$\chi \in \roman{Hom}(F,G)$
of
$r(\chi)=
[\chi(x_1),\chi(y_1)] \cdot\dots\cdot
[\chi(x_\ell),\chi(y_\ell)]$.

\proclaim{Proposition 1.8}
At a homomorphism
$\chi$ from  $F$ to $G$,
the cotangent map
$\roman T^*_{\chi}r$
from
$\roman T^*_{r(\chi)}G$
to
$\roman T^*_{\chi}\roman{Hom}(F,G)$
and the operation of left translation
make commutative the diagram
$$
\CD
\roman T^*_{\bold y}\roman{Hom}(F,G)
@<{\roman T^*_{\chi}r}<<
\roman T^*_{r(\chi)}G
\\
@V{\roman L^*_{\chi}}VV
@V{\roman L^*_{r(\chi)}}VV
\\
(g^*)^{2\ell}
@<{\de_2^{\chi}}<<
g^*
\endCD
$$
where
$\de_2^{\chi}$ refers to the corresponding operator
in
{\rm (1.7)}.
\endproclaim

\demo{Proof} This follows at once from  \cite\smooth\ (5.4). \qed
\enddemo

Recall that $IK$ denotes the augmentation ideal of a discrete group $K$.

\proclaim{Corollary 1.9}
At a homomorphism
$\phi$ from  $F$ to $G$,
having the property that $\phi(r)$
lies in the centre of $G$,
left translation
yields a commutative diagram
$$
\CD
0
@<<<
\roman{coker}\left(\roman T^*_{\phi}r\right)
@<<<
\roman T^*_{\phi}\roman{Hom}(F,G)
@<{\roman T^*_{\phi}r}<<
\roman T^*_{r(\phi)}G
\\
@.
@VVV
@V{\roman L^*_{\phi}}VV
@V{\roman L^*_{r(\phi)}}VV
\\
0
@<<<
I\pi\otimes_{\bold R \pi} g^*_{\phi}
@<<<
(g^*)^{2\ell}
@<{\de_2^{\phi}}<<
g^*
\endCD
$$
with exact columns
and hence
induces an isomorphism
from
$\roman{coker}(\roman T^*_{\phi}r)$
onto the space
$I\pi\otimes_{\bold R \pi} g^*_{\phi}$.
\endproclaim

\demo{Proof}
In fact,
the boundary operator
$\partial_1^{\phi}$
from $C_1(\Cal P,g^*_{\phi})$ to
$C_0(\Cal P,g^*_{\phi})$
induces an isomorphism
from
$
C_1(\Cal P,g^*_{\phi})\big /
\partial_2^{\phi}(C_2(\Cal P,g^*_{\phi}))
$
onto
$I\pi\otimes_{\bold R \pi} g^*_{\phi}$. \qed
\enddemo

\smallskip
Next we pass to corresponding $G$-invariant objects.

\proclaim{Lemma 1.10}
Let $f$ be a smooth real valued
$G$-invariant
function on
$\roman{Hom}(F,G)$
and
$\phi$
a  homomorphism from $F$ to $G$,
having the property that $\phi(r)$
lies in the centre of $G$.
Then
under the dual
left transformation
$\roman L^*_{\phi}$
from
$\roman T^*_{\phi}\roman{Hom}(F,G)$
to
$(g^*)^{2\ell}$
the differential
$df(\phi) \in \roman T^*_{\phi}\roman{Hom}(F,G)$
goes to
a cycle, that is,
$\roman L^*_{\phi}df(\phi)=df(\phi)\circ\roman L_{\phi}$
lies in the subspace $Z_1(\Cal P,g^*_{\phi})$
of cycles in $C_1(\Cal P,g^*_{\phi})=(g^*)^{2\ell}$
and hence determines a class
$[df(\phi)\circ\roman L_{\phi}] \in \roman H_1(\pi,g^*_{\phi})$.
\endproclaim

\demo{Proof}
The value
$\partial_1^{\phi}(df(\phi)\circ\roman L_{\phi})$ lies in
$C_0(\Cal P,g^*_{\phi})=g^*$.
Let $X \in g$ and
$$
\bold w=
(a_1,b_1,\dots,a_\ell,b_\ell)=
\left(\phi(x_1),\dots,\phi(y_\ell)\right) \in G^{2\ell}.
$$
In view of the description of (1.2) given above,
$$
\left(\partial_1^{\phi}(df(\phi)\circ\roman L_{\phi})\right)X
=
\left(df(\phi)\circ
\roman L_{\phi}\right)
\left(
X-\roman{Ad}(a^{-1}_1) X ,
\dots,
X-\roman{Ad}(b^{-1}_\ell) X
\right).
$$
However,
at
$\phi$ which, in the present description,
amounts to $\bold w$,
the analytic path
$$
t \longmapsto
\left(\roman{exp}(-tX) a_1 \roman{exp}(tX),
\dots,
\roman{exp}(-tX) b_\ell \roman{exp}(tX)\right)
$$
has
tangent vector
$$
\roman L_{\bold w}
\left(
X-\roman{Ad}(a^{-1}_1) X ,
\dots,
X-\roman{Ad}(b^{-1}_\ell) X
\right)
\in \roman T_{\bold w}G^{2\ell} \cong\roman T_{\phi}\roman{Hom}(F,G).
$$
Since $f$ is $G$-invariant,
it is constant along this path
whence its derivative
along the tangent vector
to this path is zero.
Hence
$(\partial_1^{\phi}df(\phi) \circ \roman L_{\phi})X = 0$.
Since $X \in g$ was arbitrary,
$
\partial_1^{\phi}(df(\phi)\circ \roman L_{\phi}) = 0
$
as asserted. \qed
\enddemo

\smallskip
We now
return to our principal bundle
$\xi\colon P \to \Sigma$ over $\Sigma$
with structure group $G$.
As in the predecessors
\cite\singula\ -- \cite\smooth\
to the present paper,
we choose
an invariant scalar product on $g$
and a Riemannian metric
on $\Sigma$.
We then pick
smooth closed curves
$v_1,w_1,\dots,v_\ell,w_\ell$
in $\Sigma$
representing the generators
$x_1,y_1,\dots,x_\ell,y_\ell$
so that the standard cell decomposition of
$\Sigma$ with a single 2-cell $e$
corresponding to $r$ results;
further,
let $\widehat Q \in P$
be a base point
so that $\xi (\widehat Q) =Q$.
Then the holonomy along these paths yields
the Wilson loop
mapping
$\rho$ from $\Cal A(\xi)$ to  $\roman{Hom}(F,G)$,
cf. Section 2 of~\cite\smooth.
The image
$\rho(\Cal N(\xi))$
in
$\roman{Hom}(F,G)$
of the subspace
$\Cal N(\xi)$ of central Yang-Mills connections
is a space
$\roman{Hom}_{\xi}(\Gamma,G)$
of homomorphisms from the universal central extension
$\Gamma$ of $\pi$ to $G$. See
\cite\topology\ and Section 3 of~\cite\smooth.

\smallskip
Let $\phi \in \roman{Hom}_{\xi}(\Gamma,G)$.
Recall that
the Lie bracket on $g$ induces a graded bracket
$[\cdot,\cdot]_{\phi}$
on $\roman H^*(\pi,g_{\phi})$
that endows the latter with a structure
of a graded Lie algebra.
Further,
the chosen orthogonal structure on $g$
induces a graded non-degenerate bilinear
pairing
$(\cdot,\cdot)_{\phi}$
between
$\roman H^*(\pi,g_{\phi})$
and
$\roman H^{2-*}(\pi,g_{\phi})$
which, in degree 1, amounts to a
symplectic structure
$\sigma_{\phi}$
on
$\roman H^1(\pi,g_{\phi})$;
moreover,
the assignment
to $\eta \in \roman H^1(\pi,g_{\phi})$
of
$\frac 12 [\eta,\eta]_{\phi}
\in \roman H^2(\pi,g_{\phi})$
yields a momentum mapping
$\Theta_{\phi}$ from
$\roman H^1(\pi,g_{\phi})$
to
$\roman H^2(\pi,g_{\phi})$,
for the
action
of the stabilizer
$Z_{\phi} \subseteq G$
of
$\phi \in \roman{Hom}_{\xi}(\Gamma,G)$
on
$\roman H^1(\pi,g_{\phi})$.
See Section 1 of
\cite\singula\ where this is spelled out for central Yang-Mills connections
and Section 6 of \cite\smooth.
Let
$\roman{Hom}_{\xi}(\Gamma,G)^-$
denote
the subspace
of $\roman{Hom}_{\xi}(\Gamma,G)$
consisting
of points $\phi$
so that
the
operation
$$
[\cdot,\cdot]_{\phi}
\colon
\roman H^1(\pi,g_{\phi})
\otimes
\roman H^1(\pi,g_{\phi})
@>>>
\roman H^2(\pi,g_{\phi})
\tag1.11.1
$$
is zero.
Notice that
$\roman{Hom}_{\xi}(\Gamma,G)^-$
depends on the chosen orthogonal structure on $g$.

\proclaim{Proposition 1.12}
The subspace $\roman{Hom}_{\xi}(\Gamma,G)^-$
is a smooth submanifold
of \linebreak
$\roman{Hom}(F,G)$.
Moreover,
for every
$\phi \in \roman{Hom}_{\xi}(\Gamma,G)^-$,
the tangent space
$\roman T_\phi\roman{Hom}_{\xi}(\Gamma,G)^-$
coincides with the kernel of the derivative
$dr(\phi)$ from
$\roman T_{\phi}\roman{Hom}(F,G)$
to $\roman T_{\roman{exp}(X_{\xi})}G$,
and hence the spaces of differentials constitute
an exact sequence
$$
0
@<<<
\roman T^*_{\phi}\roman{Hom}_{\xi}(\Gamma,G)^-
@<<<
\roman T^*_{\phi}\roman{Hom}(F,G)
@<{d\rel^*(\phi)}<<
\roman T^*_{\roman{exp}(X_{\xi})}G.
$$
\endproclaim

The proof will be given after that of (1.13) below.
Let
$\Cal N^-(\xi)$
be the subspace of
$\Cal N(\xi) $
consisting of central Yang-Mills connections $A$ with the property
that
the operation
$$
[\cdot,\cdot]_A
\colon
\roman H^1_A(\Sigma,\roman{ad}(\xi))
\otimes
\roman H^1_A(\Sigma,\roman{ad}(\xi))
@>>>
\roman H^2_A(\Sigma,\roman{ad}(\xi))
\tag1.11.2
$$
is zero;
by \cite\singula~(2.8),
the space $\Cal N^-(\xi)$ is a
smooth submanifold of $\Cal A(\xi)$ having at
$A \in \Cal N^-(\xi)$
tangent space
$\roman T_A \Cal N^-(\xi)$
equal to the space of 1-cocycles
$Z^1_A(\Sigma,\roman{ad}(\xi))$
in $\Omega^1(\Sigma,\roman{ad}(\xi))=\roman T_A\Cal A(\xi)$.
The  group of {\it based\/} gauge transformations
is written $\Cal G^Q(\xi)$.

\proclaim{Proposition 1.13}
The Wilson loop mapping
passes to a
smooth principal $\Cal G^Q(\xi)$-fibre bundle
$\rho^-\colon\Cal N^-(\xi) \to\roman{Hom}_{\xi}(\Gamma,G)^-$.
Consequently
$\roman{Hom}_{\xi}(\Gamma,G)^-$
is a smooth submanifold
of
$\roman{Hom}(F,G)$
in such a way that, at every
$\phi \in \roman{Hom}_{\xi}(\Gamma,G)^-$,
left translation
$\roman L_{\phi}$
from
$g^{2\ell} = C^1(\Cal P,g_{\phi})$ to
$\roman T_\phi \roman{Hom}(F,G)$
identifies the space
$Z^1(\pi,g_{\phi})$
of 1-cocycles
with the usual tangent space
$\roman T_\phi \roman{Hom}_{\xi}(\Gamma,G)^-$,
viewed as a subspace of
$\roman T_\phi \roman{Hom}(F,G)$.
Moreover
$\roman{Hom}_{\xi}(\Gamma,G)^-$
is dense in
$\roman{Hom}_{\xi}(\Gamma,G)$.
\endproclaim

\demo{Proof}
The action of
$\Cal G^Q(\xi)$
on $\Cal A(\xi)$
and hence
on
$\Cal N(\xi)$
is free,
and the Wilson loop mapping $\rho$ passes to a homeomorphism
from $\Cal N(\xi)\big/\Cal G^Q(\xi)$ onto
$\roman{Hom}_{\xi}(\Gamma,G)$.
See our papers \cite\topology\
and \cite\smooth\ for details.
Moreover, given a central Yang-Mills connection $A$,
with \linebreak
$\phi = \rho(A) \in   \roman{Hom}_{\xi}(\Gamma,G)$,
twisted integration yields an isomorphism
from
$\roman H^*_A(\Sigma,\roman{ad}(\xi)$ onto
$\roman H^*(\pi,g_{\phi})$
compatible with the relevant structure,
cf. Section 4 of~\cite\smooth.
Consequently
the
operation (1.11.1)
is zero
if and only if
(1.11.2)
is zero whence the
restriction $\rho^-$ of
the Wilson loop mapping $\rho$
is a principal
fibre bundle projection map,
 manifestly smooth, having
at every
$\phi \in \roman{Hom}_{\xi}(\Gamma,G)^-$
the asserted derivative.
Finally,
cf. our paper \cite\singulat,
the pre-image
$\Cal N^{\roman {top}}(\xi)
\subseteq
\Cal N(\xi)$
of the top stratum
$N^{\roman {top}}(\xi)$
of $N(\xi)$
is contained
in $\Cal N^-(\xi)$,
and,
by \cite\singulat~(1.4), the subspace
$N^{\roman {top}}(\xi)$
is dense in
$N(\xi)$. This implies that
$\roman{Hom}_{\xi}(\Gamma,G)^-$
is dense in
$\roman{Hom}_{\xi}(\Gamma,G)$. \qed
\enddemo

\demo{Proof of {\rm (1.12)}}
In view of \cite\smooth~(5.4),
which includes the statement dual to (1.8) above,
(1.13) implies at once the statement of (1.12). \qed
\enddemo

\proclaim{Corollary 1.14}
Let $f$ be a smooth real valued
function on
$\roman{Hom}(F,G)$ that vanishes on
$\roman{Hom}_{\xi}(\Gamma,G)^-$.
Then its differential
$
df(\phi) \in
\roman T^*_{\phi}\roman{Hom}(F,G)
$
at $\phi \in \roman{Hom}_{\xi}(\Gamma,G)^-$
passes to zero in
$
\roman T^*_{\phi}\roman{Hom}_{\xi}(\Gamma,G)^-
$.
\endproclaim

\proclaim{Corollary 1.15}
Let $f$ be a smooth real valued
$G$-invariant function on
$\roman{Hom}(F,G)$ that vanishes on
$\roman{Hom}_{\xi}(\Gamma,G)^-$.
Then,
cf. {\rm (1.10)},
the homology class
$
[df(\phi)
\circ
\roman L_{\phi}]
$
in
$\roman H_1(\pi,g^*_{\phi})$
determined by
its differential
$
df(\phi) \in
\roman T^*_{\phi}\roman{Hom}_{\xi}(\Gamma,G)^-
$
at $\phi \in \roman{Hom}_{\xi}(\Gamma,G)^-$
is zero.
\endproclaim

\smallskip
For intelligibility we recall that, by definition,
${
\roman{Rep}_{\xi}(\Gamma,G)
=
\roman{Hom}_{\xi}(\Gamma,G)\big/ G
}$
and we reproduce the construction of the smooth structure
$C^{\infty}\left(\roman{Rep}_{\xi}(\Gamma,G)\right)$;
see Section 3 of \cite\smooth\ for details.
Let
$I_{\xi}$ denote the ideal
in the algebra $C^{\infty}\left(\roman{Hom}(F,G)\right)$
of smooth functions
on $\roman{Hom}(F,G)$
that vanish on the subspace
$\roman{Hom}_{\xi}(\Gamma,G)$
of $\roman{Hom}(F,G)$.
The algebra of
{\it Whitney smooth functions\/}
$C^{\infty}\left(\roman{Hom}_{\xi}(\Gamma,G)\right)=
C^{\infty}\left(\roman{Hom}(F,G)\right)/I_{\xi}$
endows
$\roman{Hom}_{\xi}(\Gamma,G)$
with a {\it smooth structure\/}.
We then {\it define\/} the
{\it smooth structure\/}
of $\roman{Rep}_{\xi}(\Gamma,G)$
to be the algebra
$$
C^{\infty}\left(\roman{Rep}_{\xi}(\Gamma,G)\right)
= \left(C^{\infty}(\roman{Hom}(F,G))\right)^G\big/I_{\xi}^G
$$
of smooth $G$-invariant functions $\left(C^{\infty}(\roman{Hom}(F,G))\right)^G$
on $\roman{Hom}(F,G)$
modulo its ideal
$I_{\xi}^G$
of
functions
that vanish on
$\roman{Hom}_{\xi}(\Gamma,G)$.
By construction this is an algebra of functions
on
$\roman{Rep}_{\xi}(\Gamma,G)$ in an obvious fashion.
\smallskip
For $[\phi] \in \roman{Rep}_{\xi}(\Gamma,G)$,
the dual of
the real vector space
$\Omega_{[\phi]}\roman{Rep}_{\xi}(\Gamma,G)$
of differentials at $[\phi]$
with reference to
$C^{\infty}\roman{Rep}_{\xi}(\Gamma,G)$
is, by definition,
the {\it Zariski tangent space\/}
$\roman T_{[\phi]}\roman{Rep}_{\xi}(\Gamma,G)$.
At a singular
point $\phi \in \roman{Hom}_{\xi}(\Gamma,G)$,
the statements
of (1.14) and (1.15)
are still true
when the tangent space is replaced by
the Zariski tangent space
with reference to the smooth structure
$C^{\infty}\roman{Hom}_{\xi}(\Gamma,G)$.
This
fact will not be needed here
and we refrain from spelling out details;
it follows from
\cite\smooth~(7.14).
\smallskip
Let $[\phi] \in \roman{Rep}_{\xi}(\Gamma,G)$.
By \cite\smooth\ (7.10),
a choice of representative
$\chi \in \roman{Hom}_{\xi}(\Gamma,G)$
induces a linear map
$\lambda_{\chi}$ from
$\roman H^1(\pi,g_{\chi})$
to
$\roman T_{[\phi]}\roman{Rep}_{\xi}(\Gamma,G)$
which is an isomorphism if and only if
$[\phi]$ is a non-singular point of $\roman{Rep}_{\xi}(\Gamma,G)$;
the dual
$$
\lambda^*_{\chi}
\colon
\Omega_{[\phi]}\roman{Rep}_{\xi}(\Gamma,G)
=\roman T^*_{[\phi]}\roman{Rep}_{\xi}(\Gamma,G)
@>>>
\roman H_1(\pi,g^*_{\chi})
\tag1.16
$$
furnishes then a linear map
from
$\Omega_{[\phi]}\roman{Rep}_{\xi}(\Gamma,G)$
to
$\roman H_1(\pi,g^*_{\chi})$
which is as well an isomorphism if and only if
$[\phi]$ is a non-singular point of $\roman{Rep}_{\xi}(\Gamma,G)$,
and \cite\smooth\ (7.10) entails at once the following.

\proclaim{Proposition 1.17}
The linear map
$\lambda^*_{\chi}$ is independent of the choice of $\chi$ in the sense that,
given $x \in G$,
the corresponding
linear map
$\lambda^*_{x\chi}$
from
$\roman T^*_{[\phi]}\roman{Rep}_{\xi}(\Gamma,G)$
to
$\roman H_1(\pi,g^*_{x\chi})$
coincides with the composite
of $\lambda^*_{\chi}$
and the isomorphism
$\roman{Ad^*}(x)_{\flat}$
from
$\roman H_1(\pi,g^*_{\chi})$
onto
$\roman H_1(\pi,g^*_{x\chi})$ induced by $x$.
\endproclaim

\smallskip
For a  homomorphism
$\phi$ from  $F$ to $G$
having the property that $\phi(r)$
lies in the centre of $G$,
the obvious
non-degenerate real-valued bilinear pairing
between
$C^1(\Cal P,g_{\phi})$ and
$C_1(\Cal P,g^*_{\phi})$
induces a
non-degenerate bilinear pairing
between
$Z^1(\pi,g_{\phi})$ and
$I\pi\otimes_{\bold R \pi} g^*_{\phi}$
which, in turn, induces
the
{\it cap\/} pairing
$\cap$ between $\roman H^1(\pi,g_{\phi})$ and
$\roman H_1(\pi,g^*_{\phi})$;
it
is manifestly non-degenerate and
identifies
$\roman H_1(\pi,g^*_{\phi})$
with the dual of
$\roman H^1(\pi,g_{\phi})$.
This observation yields at once a proof of the following.

\proclaim{Proposition 1.18}
At a class $[\phi] \in \roman{Rep}_{\xi}(\Gamma,G)$,
for every
representative
$\chi \in [\phi]$,
the canonical evaluation pairing between
$\roman T_{[\phi]}\roman{Rep}_{\xi}(\Gamma,G)$ and
$\roman T^*_{[\phi]}\roman{Rep}_{\xi}(\Gamma,G)$
equals
the composite
$$
\roman T_{[\phi]}\roman{Rep}_{\xi}(\Gamma,G)
\otimes
\roman T^*_{[\phi]}\roman{Rep}_{\xi}(\Gamma,G)
@>{\lambda_{\chi} \otimes \lambda^*_{\chi}}>>
\roman H^1(\pi,g_{\chi})
\otimes
\roman H_1(\pi,g^*_{\chi})
@>{\cap}>>
\bold R. \qed
$$
\endproclaim

\smallskip
Let $A$ be a central Yang-Mills
connection, and let
$((\cdot,\cdot))_A$
denote the canonical evaluation pairing
between
$\roman H^1_{A}(\Sigma,\roman{ad}(\xi))$ and
$\roman H^1_{A}(\Sigma,\roman{ad}^*(\xi))$
obtained from the wedge product of forms and integration
in the usual way.
We define the
{\it dual twisted integration isomorphism\/}
$$
\roman{Int}_A^*
\colon
\roman H_1(\pi,g^*_{\rho (A)})
@>>>
\roman H^1_{A}(\Sigma,\roman{ad}^*(\xi))
\tag1.19
$$
by
$\left(\roman{Int}_A \alpha\right) \cap u
=
((\alpha , \roman{Int}_A^*u))_A$,
for $\alpha \in \roman H^1_{A}(\Sigma,\roman{ad}(\xi))$
and $u \in \roman H_1(\pi,g^*_{\rho(A)})$.
Moreover, we denote by
$\lambda^*_A$
the linear map
from
$\roman T^*_{[A]} N({\xi})$ to
$\roman H_A^1(\Sigma,\roman{ad}^*(\xi))$
which is the dual of the linear map
$\lambda_A$ given in \cite\smooth~(7.9).

\proclaim{Proposition 1.20}
For every central Yang-Mills connection $A$,
the dual twisted integration isomorphism
{\rm (1.19)}
makes commutative
the diagram
$$
\CD
\roman H_A^1(\Sigma,\roman{ad}^*(\xi))
@<{\lambda^*_A}<<
\roman T^*_{[A]} N({\xi})
\\
@A{\roman{Int}^*_A}AA
@A{d\rho_{\flat}[A]^*}AA
\\
\roman H_1(\pi,g^*_{\rho(A)})
@<{\lambda^*_{\rho(A)}}<<
\roman T^*_{[{\rho(A)}]} \roman {Rep}_{\xi}(\Gamma,G).
\endCD
$$
In particular, when $A$ represents a non-singular point,
the cotangent map
$d\rho_{\flat}[A]^*$
from
$\roman T^*_{[{\rho(A)}]} \roman {Rep}_{\xi}(\Gamma,G)$
to
$\roman T^*_{[A]} N({\xi})$
amounts
to the isomorphism
{\rm (1.19)}.
\endproclaim

\demo{Proof} This follows at once from the commutativity of
\cite\smooth~(7.11). \qed
\enddemo

\medskip\noindent{\bf 2. Poisson structures}\smallskip\noindent
Let
$
\langle \cdot,\cdot \rangle
$
be a
coadjoint action invariant
symmetric bilinear form
on $g^*$,
not necessarily non-degenerate.
For every $\phi \in \roman{Hom}_{\xi}(\Gamma,G)$,
it induces an intersection pairing
$$
\langle
\cdot,\cdot
\rangle_{\phi}
\colon
\roman H_1(\pi,g^*_{\phi})
\otimes
\roman H_1(\pi,g^*_{\phi})
\longrightarrow
\bold R.
\tag2.1.1
$$
For a smooth $G$-invariant function
$f$ on $\roman{Hom}(F,G)$,
we write
$[f] \in C^{\infty}(\roman{Rep}_{\xi}(\Gamma,G))$
for its image, obtained by restriction
of $f$ to
$\roman{Hom}_{\xi}(\Gamma,G)$.
For every $f$ and $h$
in
$\left(C^{\infty}(\roman{Hom}(F,G))\right)^G$
and every $\phi \in \roman{Hom}_{\xi}(\Gamma,G)$, let
$$
(f\bullet h) (\phi)
=
\left\langle
[df(\phi)\circ \roman L_{\phi}] \otimes [dh(\phi)\circ \roman L_{\phi}]
\right\rangle_{\phi}
\tag2.1.2
$$
where the notation
$[df(\phi)\circ \roman L_{\phi}]$ and $[dh(\phi)\circ \roman L_{\phi}]$
indicates
homology classes
in $\roman H_1(\pi,g^*_{\phi})$, cf. (1.10).
This furnishes a bilinear pairing
$$
\bullet \
\colon
\left(C^{\infty}(\roman{Hom}(F,G))\right)^G
\otimes
\left(C^{\infty}(\roman{Hom}(F,G))\right)^G
@>>>
\roman{Map}(\roman{Hom}_{\xi}(\Gamma,G),\bold R) .
\tag2.1.3
$$
By a {\it symplectic\/}
Poisson structure we mean one whose Casimir elements are the constants
only.

\proclaim{Theorem 2.1}
The pairing {\rm (2.1.3)}
induces a Poisson bracket
$$
\{\cdot,\cdot\}
\colon
C^{\infty}(\roman{Rep}_{\xi}(\Gamma,G))
\otimes
C^{\infty}(\roman{Rep}_{\xi}(\Gamma,G))
@>>>
C^{\infty}(\roman{Rep}_{\xi}(\Gamma,G))
\tag2.1.4
$$
which,
for every $f$ and $h$
in
$\left(C^{\infty}(\roman{Hom}(F,G))\right)^G$
and every $\phi \in \roman{Hom}_{\xi}(\Gamma,G)$, is calculated by the formula
$$
\{[f],[h]\}[\phi]
=
\left\langle
[df(\phi)\circ \roman L_{\phi}] \otimes [dh(\phi)\circ \roman L_{\phi}]
\right\rangle_{\phi}.
\tag2.1.5
$$
When the bilinear form
$\langle \cdot,\cdot \rangle$
is non-degenerate
(but not necessarily positive definite)
the Poisson bracket is symplectic.
\endproclaim

\smallskip\noindent
{\smc Remark.}
The description (2.1.5) of the Poisson bracket
involves choices of representatives
$f$ and $h$
of
$[f],[h] \in
C^{\infty}(\roman{Rep}_{\xi}(\Gamma,G))$, respectively,
{\it and\/} of a representative $\phi$ of the point $[\phi]$
of $\roman{Rep}_{\xi}(\Gamma,G)$.
In a sense, the choice of $\phi$ amounts to
introduction of local coordinates.
In view of the construction
of the linear map
$
\lambda_{\phi}
$
from
$\roman H^1(\pi,g_{\chi})$
to
$\roman T_{[\phi]}\roman{Rep}_{\xi}(\Gamma,G)$
in Section 7 of \cite\smooth,
it is clear
that
(2.1.5) amounts to the {\it intrinsic\/}
description (0.1) given in the Introduction.
The construction of
$
\lambda_{\phi}
$
relies on the determination of appropriate Zariski
tangent spaces
given in Section 7 of \cite\smooth\  while
our construction of the Poisson bracket,
in particular
(2.1.5) above
does {\it not\/}.
The results
in Section 7 of \cite\smooth\
are needed here {\it merely\/}
to obtain the {\it description\/}
(0.1)
of the Poisson bracket.
\smallskip
We first give an outline of the proof of (2.1).
Let $f$ and $h$
be smooth $G$-invariant functions on
$\roman{Hom}(F,G)$; we shall establish the following facts.
\newline
\noindent
(1) The function $f\bullet h$
is the restriction to
$\roman{Hom}_{\xi}(\Gamma,G)$
of a smooth $G$-invariant function
on
$\roman{Hom}(F,G)$.
\newline
\noindent
(2) When $h$ is zero on
$\roman{Hom}_{\xi}(\Gamma,G)$
so is
$f\bullet h$.
\newline
\noindent
(3)
The bracket (2.1.4) satisfies the Leibniz rule.
\newline
\noindent
(4)
The bracket (2.1.4) satisfies the Jacobi identity.
\smallskip
The proof will proceed in steps.
We shall construct a pairing
$$
\diamondsuit\
\colon
\left(C^{\infty}(\roman{Hom}(F,G))\right)
\otimes
\left(C^{\infty}(\roman{Hom}(F,G))\right)
@>>>
C^{\infty}(\roman{Hom}(F,G))
$$
satisfying the Leibniz rule;
its construction will rely
on a combinatorial description
of the intersection pairing
to be given
in (2.3) -- (2.9)
and will be completed in (2.12)
while its $G$-invariance will be established
in (2.14).
Statement (2.15) below will imply that
the resulting pairing
on $\left(C^{\infty}(\roman{Hom}(F,G))\right)^G$,
combined with the projection
from
$(C^{\infty}(\roman{Hom}(F,G)))^G$
onto
$C^{\infty}(\roman{Rep}_{\xi}(\Gamma,G))$,
comes down to (2.1.3).
In (2.17) we shall then show
that
$\diamondsuit$
passes to a pairing
on
$C^{\infty}(\roman{Rep}_{\xi}(\Gamma,G))$;
this takes care of (2) in the above outline.
The Jacobi identity
and, moreover, the
symplecticity
of the Poisson bracket
for a non-degenerate  bilinear form
on $g^*$
will be established in the next Section.
\smallskip
We now begin working out the details.
For a while we shall admit
an arbitrary commutative ring $R$
as ground ring.
We proceed
at first
towards a description of
the requisite intersection pairings.
With reference to (1.4) above, let
$$
\widetilde {\roman R}(\Cal P)
=
\roman {Hom}_{R\pi}({\roman R}(\Cal P),R\pi).
\tag2.2
$$
With the notation
$
\widetilde {\roman R}_{2-j}(\Cal P)=
\roman {Hom}_{R\pi}({\roman R}_j(\Cal P),R\pi)$,
for $0 \leq j \leq 2$, it
looks like
$$
\widetilde {\roman R}(\Cal P)
\colon
\widetilde {\roman R}_2(\Cal P)
@>{\partial_2}>>
\widetilde {\roman R}_1(\Cal P)
@>{\partial_1}>>
\widetilde {\roman R}_0(\Cal P),
\tag2.3.1
$$
and, in the standard way,
the
${\roman R}_j(\Cal P)$
coming as
{\it right\/} $R\pi$-modules,
each
$\widetilde {\roman R}_j(\Cal P)$
inherits a structure of a {\it left\/}
$R\pi$-module
by means of
$$
(x\phi)(y) = x(\phi y) ,\quad x \in R\pi, \ y \in {\roman R}_j(\Cal P),
\ j = 0,1,2 .
\tag2.3.2
$$
Furthermore
the canonical map
from
$\roman R(\Cal P)$
to
$\roman {Hom}_{R\pi}(\widetilde{\roman R}(\Cal P),R\pi)$
is an isomorphism of free
resolutions of $R$ in the category of right $R\pi$-modules
and, for every left
$R\pi$-module
$W$,
the assignment to
$f \otimes w
\in
\roman {Hom}_{R\pi}(\widetilde{\roman R}(\Cal P),R\pi)\otimes_{R\pi} W$
of
$f_w$
given by the formula
$f_w (y) = f(y) w$,
for
$y \in  \widetilde{\roman R}(\Cal P)$,
yields
an isomorphism
of chain complexes
from
$\roman {Hom}_{R\pi}(\widetilde{\roman R}(\Cal P),R\pi)\otimes_{R\pi} W$
onto
$\roman {Hom}_{R\pi}(\widetilde{\roman R}(\Cal P),W)$;
hence the resulting composite
$$
\roman R(\Cal P)\otimes_{R\pi} W
@>>>
\left(\roman {Hom}_{R\pi}(\widetilde{\roman R}(\Cal P),R\pi)\right)
\otimes_{R\pi} W
@>>>
\roman {Hom}_{R\pi}(\widetilde{\roman R}(\Cal P),W)
\tag2.3.3
$$
is likewise an isomorphism of chain complexes.
Poincar\'e duality may now be expressed in the following form.

\proclaim{Proposition 2.3}
The chain complex
$
\widetilde {\roman R}(\Cal P)
$
is a free resolution of the ground ring $R$ in the category of
left $R\pi$-modules
in such a way that,
for every left $R\pi$-module $W$,
when $\roman H^*(\pi,W)$ is calculated as
$\roman H^*(\roman{Hom}_{R\pi}(\widetilde {\roman R}(\Cal P),W))$
and
$\roman H_*(\pi,W)$ as
$\roman H_*({\roman R}(\Cal P)\otimes _{R\pi} W)$,
the Poincar\'e duality isomorphism
$$
[\pi]\cap - \ \colon
\roman H^*(\pi,W)
\longrightarrow
\roman H_{2-*}(\pi,W),\quad 0 \leq j \leq 2,
$$
that is, the cap product with the fundamental class
$[\pi]\in \roman H_2(\pi,R)$,
is induced by the
inverse of the
canonical isomorphism
{\rm (2.3.3)}.
\endproclaim

\demo{Proof}
By construction,
the homology of
$\widetilde {\roman R}(\Cal P)$
is the cohomology
of $\pi$ with values in
$R\pi$.
However,
this cohomology
is just the cohomology
$\roman H^*_f(\widetilde \Sigma,R)$
with finite cochains
of the universal covering
$\widetilde \Sigma$
and, by Poincar\'e duality,
$
\roman H^*_f(\widetilde \Sigma,R)$
is isomorphic to
$\roman H_{2-*}(\widetilde\Sigma,R)$
whence
$
\roman H^0(\pi,R\pi) = 0,
\
\roman H^1(\pi,R\pi) = 0,
\
\roman H^2(\pi,R\pi) \cong R$.
Thus
$
\widetilde {\roman R}(\Cal P)
$
is a free resolution of the ground ring $R$ in the category of
left $R\pi$-modules.
\smallskip
To get our hands on Poincar\'e duality,
let $\bold P_1,\, \bold P_2,\, \bold P_3$ be three free resolutions of $R$
in the category of left $R\pi$-modules, let
${
\Delta
\colon
\bold P_1
@>>>
\bold P_2 \otimes \bold P_3
}$
be a diagonal map,
and let  $W$ be a left
$R\pi$-module
and
$V$ a right
$R\pi$-module.
For every $\phi \in \roman{Hom}_{R\pi}(\bold P_2,W)$,
we then have
the chain map
$$
- \cap \phi
\colon
V\otimes_{R\pi} \bold P_1
@>>>
V\otimes_{R\pi} \left(W\otimes \bold P_3\right),
$$
defined as the composite
$$
V\otimes_{R\pi} \bold P_1
@>{\roman{Id} \otimes \Delta}>>
V\otimes_{R\pi} \left(\bold P_2 \otimes \bold P_3\right)
@>{\roman{Id} \otimes \phi \otimes \roman{Id}}>>
V\otimes_{R\pi} \left(W\otimes \bold P_3\right) ;
$$
upon taking $V = R$ and writing
$\overline {\bold P_1} =R\otimes_{R\pi} \bold P_1$,
we arrive at a chain map
$$
\overline {\bold P_1}
@>>>
\roman{Hom}\left(
\roman{Hom}_{R\pi}(\bold P_2,W),
R\otimes_{R\pi} \left(W\otimes \bold P_3\right)
\right)
$$
which assigns
$a \cap -$ to
$a \in \overline {\bold P_1}$.
Moreover,
write
$\bold P^{\roman r}_3$
for the free resolution
of $R$ in the category of right $R\pi$-modules
obtained from
$\bold P_3$ in the usual way, that is,
as a chain complex in the category of $R$-modules,
$\bold P^{\roman r}_3 = \bold P_3$,
and the group $\pi$ acts on the right
in the obvious way
so that
$y x = x^{-1}y$
for every $x\in \pi$ and every  $y \in \bold P_3$;
then
the chain map
$$
R\otimes_{R\pi} \left(W\otimes \bold P_3\right)
@>>>
\bold P^{\roman r}_3 \otimes_{R\pi} W
$$
given by
${
1\otimes w \otimes y \mapsto y \otimes w,
\quad
w \in W,\ y \in \bold P_3,
}$
is an isomorphism.
Consequently
the assignment
to $a \in \overline {\bold P_1}$
of $a \cap -$
yields a chain map
$$
\overline {\bold P_1}
@>>>
\roman{Hom}\left(
\roman{Hom}_{R\pi}(\bold P_2,W),
\bold P^{\roman r}_3 \otimes_{R\pi} W
\right)
$$
which, by construction, induces the operation
$$
\roman H_i(\pi,R)
@>>>
\roman{Hom}\left(
\roman H^*(\pi,W),
\roman H_{i-*}(\pi,W)
\right)
$$
sending
$a\in \roman H_i(\pi,R)$ to  $a \cap -$;
in particular, with $a = [\pi] \in \roman H_2(\pi, R)$,
the fundamental class,
we get the Poincar\'e
duality isomorphism
from
$\roman H^*(\pi,W)$ onto
$\roman H_{2-*}(\pi,W)$.
\smallskip
Let $\bold P_1$ be an arbitrary
free resolution
of $R$ in the category of {\it left\/}
$R\pi$-modules,
$\bold P_2 = \widetilde {\roman R}(\Cal P)$
cf. (2.2) above,
and
$\bold P_3 = \roman R(\Cal P)^{\roman l}$,
the free resolution
$\roman R(\Cal P)$,
converted into one in the category of
{\it left\/} $R\pi$-modules
by the same kind of construction
as that used above for the passage from
left to right modules.
Then
$R\otimes_{R\pi} \left(W\otimes \bold P_3\right)$
looks like
${\roman R}(\Cal P)
\otimes_{R\pi} W$; further
$\roman{Hom}_{R\pi}(\bold P_2,W)
=
\roman{Hom}_{R\pi}(\widetilde {\roman R}(\Cal P),W)$
amounts to
${\roman R}(\Cal P)
\otimes_{R\pi} W$;
and the cap product with the fundamental class
boils down to the
inverse of the
canonical isomorphism
{\rm (2.3.3)} as asserted. \qed
\enddemo
\smallskip
Next we consider the lifted object
${
\widehat{\widetilde {\roman R}(\Cal P)}
=
\roman {Hom}_{RF}(\widehat{{\roman R}(\Cal P)},RF).
}$
With the notation
$
\widehat{\widetilde {\roman R}_{2-j}(\Cal P)}=
\roman {Hom}_{RF}(\widehat{{\roman R}_j(\Cal P)},RF),
$
it looks like
$$
\widehat{\widetilde {\roman R}(\Cal P)}
\colon
\widehat{\widetilde {\roman R}_2(\Cal P)}
@>{\widehat{\partial}_2}>>
\widehat{\widetilde {\roman R}_1(\Cal P)}
@>{\widehat{\partial}_1}>>
\widehat{\widetilde {\roman R}_0(\Cal P)},
\tag2.4
$$
and each
$\widehat{\widetilde {\roman R}_j(\Cal P)}$
inherits a structure of a left
$RF$-module
by means of a formula
of the kind
(2.3.2).
Let
${
\Delta
\colon
\widetilde {\roman R}(\Cal P)
@>>>
\widetilde {\roman R}(\Cal P)
\otimes
\widetilde {\roman R}(\Cal P)
}$
be an
$R\pi$-linear
diagonal map
for the
free resolution $\widetilde {\roman R}(\Cal P)$;
as usual, the group ring $R\pi$
acts here on the tensor product
$\widetilde {\roman R}(\Cal P)
\otimes
\widetilde {\roman R}(\Cal P)$
through the diagonal map
$\Delta \colon \pi \to \pi \times \pi$.
We lift the diagonal map $\Delta$ to an
$RF$-linear
morphism
${
\widehat {\Delta}
\colon
\widehat{\widetilde {\roman R}(\Cal P)}
@>>>
\widehat{\widetilde {\roman R}(\Cal P)}
\otimes
\widehat{\widetilde {\roman R}(\Cal P)}
}$
of graded modules so that the diagram
$$
\CD
\widehat{\widetilde {\roman R}(\Cal P)}
@>{\widehat \Delta}>>
\widehat{\widetilde {\roman R}(\Cal P)}
\otimes
\widehat{\widetilde {\roman R}(\Cal P)}
\\
@VVV
@VVV
\\
\widetilde {\roman R}(\Cal P)
@>{\Delta}>>
\widetilde {\roman R}(\Cal P)
\otimes
\widetilde {\roman R}(\Cal P)
\endCD
$$
is commutative.
Thus $\widehat {\Delta}$ is a kind of diagonal map for
$\widehat{\widetilde {\roman R}(\Cal P)}$.
\smallskip
Let $U$ be a left $RF$-module,
with structure map
$\chi$ from  $F$ to $\roman{Aut}(U)$; then
the component
${
\widehat \Delta
\colon
\widehat{\widetilde {\roman R}_2(\Cal P)}
@>>>
\widehat{\widetilde {\roman R}_1(\Cal P)}
\otimes
\widehat{\widetilde {\roman R}_1(\Cal P)}
}$
induces
a pairing
$$
m_{\chi}
\colon
\roman{Hom}_{RF}(\widehat{\widetilde {\roman R}_1(\Cal P)},U)
\otimes
\roman{Hom}_{RF}(\widehat{\widetilde {\roman R}_1(\Cal P)},U)
@>>>
\roman{Hom}_{RF}(\widehat{\widetilde {\roman R}_2(\Cal P)},U\otimes U),
\tag2.5
$$
the tensor product
$U\otimes U$ being equipped with the diagonal
$RF$-module structure.
Moreover,
given another left
$RF$-module $V$,
with structure map
$\theta$ from  $F$ to $\roman{Aut}(V)$,
and a morphism
$\alpha \colon U \to V$ of left
$RF$-modules,
by naturality,
the resulting maps
$m_{\chi}$ and
$m_{\theta}$ are compatible in the sense that the corresponding diagram
$$
\CD
\roman{Hom}_{RF}(\widehat{\widetilde {\roman R}_1(\Cal P)},U)
\otimes
\roman{Hom}_{RF}(\widehat{\widetilde {\roman R}_1(\Cal P)},U)
@>{m_{\chi}}>>
\roman{Hom}_{RF}(\widehat{\widetilde {\roman R}_2(\Cal P)},U\otimes U)
\\
@VVV
@VVV
\\
\roman{Hom}_{RF}(\widehat{\widetilde {\roman R}_1(\Cal P)},V)
\otimes
\roman{Hom}_{RF}(\widehat{\widetilde {\roman R}_1(\Cal P)},V)
@>{m_{\theta}}>>
\roman{Hom}_{RF}(\widehat{\widetilde {\roman R}_2(\Cal P)},V\otimes V),
\endCD
\tag2.6
$$
is commutative.
\smallskip
Taking into account the above isomorphism
(2.3.3)
we see that,
$U^{2\ell}$ being identified with
${\roman R}_1(\Cal P)\otimes _{R\pi} U$
and
${\roman R}_0(\Cal P)\otimes _{R\pi}(U\otimes U)$
with $U\otimes U$,
the pairing (2.5)
looks like
$$
m_{\chi}
\colon
U^{2\ell} \otimes U^{2\ell}
@>>>
U\otimes U.
\tag2.7
$$
By construction,
when the left
$RF$-module structure $\chi$ on $U$
comes from a left
$R\pi$-module structure,
the pairing (2.5)
induces the cup pairing
$$
\cup\colon
\roman H^1(\pi,U) \otimes
\roman H^1(\pi,U)
@>>>
\roman H^2(\pi,U\otimes U) ,
\tag2.8
$$
computed from the free resolution
$
\widetilde {\roman R}(\Cal P)
$
of $R$ in the category of {\it left\/} $R\pi$-modules,
and hence {\rm (2.7)}
induces
the intersection
pairing
$$
\roman H_1(\pi,U) \otimes
\roman H_1(\pi,U)
@>{\iota}>>
\roman H_0(\pi,U\otimes U) ,
\tag2.9
$$
computed from the
{\it standard\/}
free resolution  (1.4)
of $R$ in the category of right
$R\pi$-modules.
In fact, with reference to
the
left
$\pi$-module structure
$\chi \colon \pi \to \roman{Aut}(U)$
on $U$
and the corresponding one
$\chi^{\otimes} \colon \pi \to \roman{Aut}(U\otimes U)$
on the tensor product $U \otimes U$, the diagram
$$
\CD
\roman H^1(\pi,U) \otimes
\roman H^1(\pi,U)
@>{\roman{\cup}}>>
\roman H^2(\pi,U\otimes U)
\\
@V{([\pi]\cap -)\otimes ([\pi]\cap -)}VV
@VV{[\pi]\cap -}V
\\
\roman H_1(\pi,U) \otimes
\roman H_1(\pi,U)
@>{\iota}>>
\roman H_0(\pi,U\otimes U)
\endCD
$$
is commutative.

\smallskip
We now take $R = \bold R$, the reals.
Given
a homomorphism $\chi$ from $F$ to $G$ as before
and taking $U = g^*_{\chi}$,
with respect to the canonical identifications
of $\roman{Hom}_{RF}
(\widehat{\widetilde {\roman R}_1(\Cal P)},g^*_{\chi})$
with
$(g^*)^{2\ell}$
and of
$\roman{Hom}_{RF}(\widehat{\widetilde {\roman R}_2(\Cal P)},
g^*_{\chi}\otimes g^*_{\chi})$
with
$g^* \otimes g^*$,
the resulting pairing (2.7) looks like
${
m_{\chi}
\colon
\left(g^*\right)^{2\ell}
\otimes
\left(g^*\right)^{2\ell}
\longrightarrow
g^* \otimes g^*
}$;
it is clear that,
on each connected component,
the resulting map
$$
m \colon
\roman{Hom}(F,G)
\longrightarrow
\roman{Hom}
\left(
\left(g^*\right)^{2\ell}\otimes \left(g^*\right)^{2\ell},
g^* \otimes g^*
\right)
\tag2.10
$$
which sends $\chi \in \roman{Hom}(F,G)$ to  $m_{\chi}$
is algebraic and hence smooth.
\smallskip
In the usual way, the group $G$ acts on
$\roman{Hom}(F,G)$
and
$\roman{Hom}((g^*)^{2\ell}\otimes (g^*)^{2\ell},
g^* \otimes g^*)$;
for an element $\beta$ of the latter, for $x \in G$, and
for $a,b \in (g^*)^{2\ell}$,
$$
(x \beta)(a \otimes b)
=
x \left(\beta(x^{-1}a \otimes x^{-1}b)\right).
$$

\proclaim{2.11}
The map $m$
is $G$-equivariant.
\endproclaim

\demo{Proof}
Given $x \in G$,
whatever $\chi \in \roman{Hom}(F,G)$,
the induced linear map
$
\roman {Ad}^*(x)
$
from
$g_{\chi}^*$ to
$g_{x\chi}^*$
is an isomorphism of $\bold RF$-modules.
By naturality,
the diagram
$$
\CD
\roman T^*_{\chi}\roman{Hom}(F,G)
\otimes
\roman T^*_{\chi}\roman{Hom}(F,G)
@>{m_\chi}>>
g_{\chi}^*
\otimes
g_{\chi}^*
\\
@V{\roman {Ad}^*(x)\otimes \roman {Ad}^*(x)}VV
@V{\roman {Ad}^*(x)\otimes \roman {Ad}^*(x)}VV
\\
\roman T^*_{x\chi}\roman{Hom}(F,G)
\otimes
\roman T^*_{x\chi}\roman{Hom}(F,G)
@>{m_{x\chi}}>>
g_{x\chi}^*
\otimes
g_{x\chi}^*
\endCD
$$
is therefore commutative.
This implies the assertion. \qed
\enddemo

\smallskip
For
arbitrary smooth functions
$f$ and $h$
on $\roman{Hom}(F,G)$
and every $\chi \in \roman{Hom}(F,G)$,
let
$$
(f\diamondsuit h)
(\chi)
=
\left\langle m_\chi\left(
(df(\chi) \circ \roman L_{\chi})
\otimes
(dh(\chi)\circ \roman L_{\chi})
\right)\right\rangle.
\tag2.12
$$
This yields
a bilinear pairing
$\diamondsuit$
on $C^{\infty}(\roman{Hom}(F,G))$
with values in
$C^{\infty}(\roman{Hom}(F,G))$.
We list some of its properties.

\proclaim{2.13}
It satisfies the Leibniz rule
${
f \diamondsuit (hk)
=
h (f\diamondsuit k) +
k(f\diamondsuit h)  ,
}$
whatever
smooth functions
$f,h,k$ on $\roman{Hom}(F,G)$.
\endproclaim

This follows at once from the construction
of $\diamondsuit$ in terms of differentials.
For
smooth $G$-invariant functions
$f$ and $h$
on $\roman{Hom}(F,G)$, we now spell out
(2.14), (2.15), and (2.17) below.

\proclaim{2.14}
The function
$f\diamondsuit h $ is
also $G$-invariant.
\endproclaim

\demo{Proof}
Let $x \in G$
and
$\chi \in \roman{Hom}(F,G)$.
Because
$f$ and $h$
are $G$-invariant real valued functions,
$
df(x\chi) = \roman{Ad}^*(x)df(\chi)
$
and
$dh(x\chi) = \roman{Ad}^*(x)dh(\chi)$.
In view of (2.11),
$$
\aligned
(f\diamondsuit h)
(x\chi)
&=
\left\langle m_{x\chi}
\left(
(df(x\chi)\circ \roman L_{x\chi})
 \otimes
(dh(x\chi)\circ \roman L_{x\chi})
\right)\right\rangle
\\
&=
\left\langle m_{x\chi}
\left(
\roman{Ad}^*(x)
(df(\chi) \circ \roman L_{\chi})
\otimes
\roman{Ad}^*(x)
(dh(\chi) \circ \roman L_{\chi})
\right)\right\rangle
\\
&=
\left\langle
(\roman{Ad}^*(x) \otimes \roman{Ad}^*(x))
m_\chi\left(
(df(\chi)\circ \roman L_{\chi})
 \otimes
(dh(\chi)\circ \roman L_{\chi})
\right)\right\rangle
\\
&=
\left\langle m_\chi
\left(
(df(\chi)\circ \roman L_{\chi})
 \otimes
(dh(\chi)\circ \roman L_{\chi})
\right)\right\rangle
\\
&=
(f\diamondsuit h)
(\chi),
\endaligned
$$
where we exploited the fact that the
given
bilinear form
$\langle \cdot,\cdot\rangle$ on $g^*$
is coadjoint action invariant. \qed
\enddemo

Recall from  {\rm (1.10)} that
$
df(\chi)\circ \roman L_{\chi} \in Z_1(\Cal P,g_{\chi}^*)
$
and
$dh(\chi) \circ \roman L_{\chi}\in Z_1(\Cal P,g_{\chi}^*)$,
whatever
$\chi \in \roman{Hom}(F,G)$.

\proclaim{2.15}
For every
$\phi \in \roman{Hom}_{\xi}(\Gamma,G)$,
the value
$(f\diamondsuit h)  (\phi)$
equals the right-hand side
${
\langle [df(\phi)\circ \roman L_{\phi}]
\otimes
[dh(\phi)\circ \roman L_{\phi}]
\rangle_{\phi}
}$
of {\rm (2.1.2)}
and hence depends only on
the classes
\linebreak
$[df(\phi)\circ \roman L_{\phi}]\in \roman H_1(\pi,g_{\phi}^*)$
and
$[dh(\phi)\circ \roman L_{\phi}]\in \roman H_1(\pi,g_{\phi}^*)$.
Consequently the pairing $\diamondsuit$,
combined with the projection
from
$(C^{\infty}(\roman{Hom}(F,G)))^G$
onto
$C^{\infty}(\roman{Rep}_{\xi}(\Gamma,G))$,
comes down to the pairing {\rm (2.1.3)}.
\endproclaim

This is a consequence of the following.

\proclaim{2.16}
For an arbitrary left $\pi$-module
$U$,
with structure map
$\chi$ from  $\pi$ to $\roman{Aut}(U)$,
the restriction of the composite
$U^{2\ell} \otimes U^{2\ell}\to \bold R$
of {\rm (2.7)} with a $\pi$-invariant bilinear form on $U$
to the subspace
${Z_1({\roman R}(\Cal P),U)
\otimes
Z_1({\roman R}(\Cal P),U)
}$
of
${({\roman R}_1(\Cal P)\otimes _{R\pi} U)\otimes
({\roman R}_1(\Cal P)\otimes _{R\pi} U)}$
factors through
the corresponding
$\bold R$-valued intersection pairing on
$\roman H_1\left(\pi,U\right)\otimes\roman H_1\left(\pi,U\right)$.
\endproclaim

\demo{Proof}
This is just another way of
saying that, with reference to the
free resolution
$\widetilde{\roman R}(\Cal P)$,
the pairing induced by the diagonal map,
when restricted to the cocycles,
factors through cohomology. \qed
\enddemo

\proclaim{2.17}
If $h$ is zero on
$\roman{Hom}_{\xi}(\Gamma,G)$,
so is the function
$f\diamondsuit h $.
In other words,
for general $h$,
the set function
$\{f,h\}$
defined on
$\roman{Rep}_{\xi}(\Gamma,G)$
by the bracket  {\rm (2.1.4)}
is well defined on
$C^{\infty}(\roman{Rep}_{\xi}(\Gamma,G))$.
\endproclaim

\demo{Proof}
Since $h$ is zero on
$\roman{Hom}_{\xi}(\Gamma,G)$,
at
$\phi \in\roman{Hom}_{\xi}^-(\Gamma,G)$,
by virtue of  (1.14),
$[dh(\phi)\circ \roman L_{\phi}] \in \roman H_1(\pi,g_{\phi}^*)$
is zero.
Consequently the value
$$
(f\diamondsuit h ) (\phi) =
\langle
[df(\phi)\circ \roman L_{\phi}]
\otimes
[dh(\phi)\circ \roman L_{\phi}]
\rangle_{\phi}
$$
is zero for every
$\phi\in \roman{Hom}_{\xi}^-(\Gamma,G)$.
Since
$\roman{Hom}_{\xi}^-(\Gamma,G)$
is dense in
$\roman{Hom}_{\xi}(\Gamma,G)$,
cf. (1.4) above,
the function
$f\diamondsuit h  $
is zero on
$\roman{Hom}_{\xi}(\Gamma,G)$. \qed
\enddemo

We can now prove compatibility with
the action of the mapping class group since this does not
rely on the Jacobi identity.

\demo{Proof of Theorem 2}
Let $\phi$ be an element of
$\roman{Hom}_\xi(\Gamma,G)$ and
consider
an automorphism
$\beta$ of $\Gamma$.
The composite $\phi \beta$
is a homomorphism from
$\Gamma$ to $G$ which lies in
$\roman{Hom}_{\beta^*\xi}(\Gamma,G)$,
where
$\beta^*\xi$ equals $\xi$ (topologically) if $\beta$
preserves the orientation and equals the
(topologically)
\lq\lq opposite\rq\rq\  bundle
otherwise
(which may coincide with $\xi$).
The automorphism  $\beta$ induces a commutative diagram
$$
\CD
\Omega_{[\phi]}\roman{Rep}_{\xi}(\Gamma,G)
@>{\lambda^*_{\phi}}>>
\roman H_1(\pi,g^*_{\phi})
\\
@A{\beta_{\sharp}}AA
@A{\beta_{\sharp}}AA
\\
\Omega_{[\phi \beta]}\roman{Rep}_{\beta^*\xi}(\Gamma,G)
@>{\lambda^*_{\phi \beta}}>>
\roman H_1(\pi,g^*_{\phi\beta }).
\endCD
$$
Comparison with the formula (0.1)
shows that the induced diffeomorphism
$\beta^\sharp$ from
$\roman{Rep}_{\beta^*\xi}(\Gamma,G)$
to
$\roman{Rep}_{\xi}(\Gamma,G)$
is compatible with the brackets (2.1.4),
taken on both
$\roman{Rep}_{\beta^*\xi}(\Gamma,G)$
and
$\roman{Rep}_{\xi}(\Gamma,G)$
if necessary.
Thus the induced action of the
group of orientation
preserving outer automorphisms of $\Gamma$
on $\roman{Rep}_{\xi}(\Gamma,G)$
preserves the bracket.
However the group of outer automorphisms of $\Gamma$
is isomorphic to
the mapping class group
of $\Sigma$
in such a way that the action corresponds to
that of the mapping class group
on $\roman{Rep}(\Gamma,G)$.
The latter, restricted to orientation preserving
mapping classes, preserves each subspace
of the kined
$\roman{Rep}_{\xi}(\Gamma,G)$.
\enddemo

\medskip\noindent{\bf 3. The Jacobi identity}
\smallskip\noindent
Let
$\langle\cdot,\cdot\rangle_g$
be an   orthogonal structure
on $g$, that is,
an adjoint action invariant
scalar product.
After a choice of Riemannian metric and orientation on
$\Sigma$,
the theory
established
in \cite\atibottw\
and
in our earlier papers
\cite\singula\
--\ \cite\smooth\
is available.
In particular,
the Wilson loop mapping
$\rho_{\flat}$
from
$N(\xi)$
to
$\roman{Rep}(\Gamma,G)$
restricts to a diffeomorphism
from
$N^{\roman{top}}(\xi)$
onto
$\roman{Rep}^{\roman{top}}_{\xi}(\Gamma,G)$
whose derivative at a point
$[A] \in N^{\roman{top}}(\xi)$
amounts to the twisted integration mapping
$\roman{Int}_A$ from
$\roman H^1_A(\Sigma,\roman{ad}(\xi))$ to
$\roman H^1(\pi,g_{\phi})$;
see e.~g.
\cite\smooth\ (7.11)
for details.
Moreover the data determine a symplectic structure
on $N^{\roman{top}}(\xi)$
which,
for $[A] \in N^{\roman{top}}(\xi)$,
on the tangent space amounts to the symplectic structure
$\sigma_A$
on
$\roman H^1_A(\Sigma,\roman{ad}(\xi))$
induced by the data, cf.
\cite\atibottw\
and our earlier paper
\cite\singula;
with $\phi = \rho(A)$,
under the twisted integration mapping, this structure then
passes to
the symplectic structure
$\sigma_{\phi}$
on
$\roman H^1(\pi,g_{\phi})$
mentioned in Section 1.
In this way, the
orthogonal structure on $g$
gives rise to a symplectic structure
on
$\roman{Rep}^{\roman{top}}_{\xi}(\Gamma,G)$.
\smallskip
Henceforth we  denote
the given
coadjoint action invariant
symmetric bilinear form
on $g^*$
by
$\langle\cdot,\cdot\rangle_{g^*}$.
It induces
a 2-tensor
$\omega$
on
$N^{\roman{top}}(\xi)$.
In fact, let
$A$ be a central Yang-Mills connection
representing a point of the top stratum $N^{\roman{top}}(\xi)$,
and let $\phi = \rho(A) \in \roman{Hom}_\xi(\Gamma,G)$.
By means of the isomorphism
$\lambda_A^*$ from
$\roman T^*_{[A]}N(\xi)$
onto $\roman H^1_A(\Sigma,\roman{ad}^*(\xi))$,
cf. Section 1 above,
at $[A]$,
the tensor $\omega$
then amounts to
the 2-form
$\omega_A$ on
$\roman H^1_{A}(\Sigma,\roman{ad}^*(\xi))$
induced by
$\langle\cdot,\cdot\rangle_{g^*}$
via the wedge product of forms and integration.
Further, by (1.20),
the
cotangent map
of
the diffeomorphism
from
$N^{\roman{top}}(\xi)$
onto
$\roman{Rep}^{\roman{top}}_{\xi}(\Gamma,G)$
boils down to
the dual twisted integration mapping (1.19);
moreover
the composite
$$
\roman H_1(\pi,g^*_{\phi})
\otimes
\roman H_1(\pi,g^*_{\phi})
@>{\roman{Int}_A^*\otimes \roman{Int}_A^*}>>
\roman H^1_{A}(\Sigma,\roman{ad}^*(\xi))
\otimes
\roman H^1_{A}(\Sigma,\roman{ad}^*(\xi))
@>{\omega_A}>>
\bold R
$$
coincides with the intersection pairing
(2.1.1) induced by the given 2-form on $g^*$.
On the other hand,
via the isomorphism $\lambda_\phi^*$ given in (1.16),
the latter pairing amounts just to the 2-tensor
at $[\phi] \in\roman{Rep}^{\roman{top}}_{\xi}(\Gamma,G)$
which
corresponds to
the bracket (2.1.3)
induced by $\langle\cdot,\cdot\rangle_{g^*}$.
\smallskip
When the bilinear form  $\langle\cdot,\cdot\rangle_{g^*}$
is positive definite,
it induces an isomorphism
between $g$ and $g^*$
and hence an orthogonal structure
$
\langle\cdot,\cdot\rangle_{g}
$
on $g$ so that under this isomorphism
the two 2-forms correspond,
and the above remarks apply;
moreover
the 2-forms
$\sigma_A$ and
$\omega_A$
then correspond to each other under adjointness,
and the resulting bracket is just the corresponding
symplectic Poisson structure.
Since $A$ is arbitrary, this shows that then
the bracket
{\rm (2.1.3)}
coincides with
the symplectic Poisson
bracket
induced by the symplectic structure
on
$\roman{Rep}^{\roman{top}}_{\xi}(\Gamma,G)$
determined by
the orthogonal structure
$
\langle\cdot,\cdot\rangle_{g}
$
on $g$
and hence in particular
satisfies the Jacobi identity.
By \cite\singulat~(1.4),
the top stratum $N^{\roman{top}}(\xi)$
is dense in $N(\xi)$, and hence
$\roman{Rep}^{\roman{top}}_{\xi}(\Gamma,G)$
is dense in
$\roman{Rep}_{\xi}(\Gamma,G)$;
consequently
for a positive definite 2-form on $g^*$
the
bracket
(2.1.3)
satisfies the
Jacobi identity
everywhere
on
$\roman{Rep}_{\xi}(\Gamma,G)$.
This establishes Theorem 2.1
in this special case.
\smallskip
To handle the case of a general 2-form $\langle\cdot,\cdot\rangle_{g^*}$
we recall that,
by structure
theory, the Lie algebra $g$ decomposes uniquely
into a direct sum
of its centre $z$
and
the  simple ideals in the semi simple Lie algebra
$[g,g]$.
This implies that
$g^*$ decomposes as a direct sum
$$
g^* = g_+^*\oplus g^*_- \oplus g^*_0
\tag3.1
$$
of $G$-modules,
together with
coadjoint action invariant
scalar products
$\langle\cdot,\cdot\rangle^*_+$ and $\langle\cdot,\cdot\rangle^*_-$
on $g_+^*$
and
$g_-^*$,
respectively,
so that  the 2-form on $g^*$
decomposes as
$$
\langle\cdot,\cdot\rangle_{g^*}
=
\langle\cdot,\cdot\rangle^*_+
-
\langle\cdot,\cdot\rangle^*_- ,
$$
and so that
$
g_0^*
$
is its
null space. We
note that, even when $G$ is not connected, the decomposition
(3.1) is one of $G$-modules; this relies on the uniqueness
of the decomposition of $g$
since $\langle\cdot,\cdot\rangle_{g^*}$ is assumed $G$-invariant.
Picking a coadjoint action invariant
scalar product
$
\langle\cdot,\cdot\rangle^*_0
$
on $g_0^*$,
we obtain a  coadjoint action invariant scalar product
$$
\langle\cdot,\cdot\rangle'_{g^*}
=
\langle\cdot,\cdot\rangle^*_+
+
\langle\cdot,\cdot\rangle^*_-
+
\langle\cdot,\cdot\rangle^*_0
\colon
g^*
\otimes
g^*
@>>>
\bold R
$$
on $g^*$;
it induces an isomorphism
from $g^*$ onto $g$
which, in turn,
identifies
$\langle\cdot,\cdot\rangle'_{g^*}$
with
an orthogonal structure
$$
\langle\cdot,\cdot\rangle'_{g}
\colon
g
\otimes
g
@>>>
\bold R
\tag3.2
$$
on $g$.
By construction, then,
the
direct sum decomposition
(3.1)
passes to a
direct sum decomposition
$$
g = g_+\oplus g_- \oplus g_0
\tag3.3
$$
of $G$-modules;
the forms
$\langle\cdot,\cdot\rangle^*_+$ and $\langle\cdot,\cdot\rangle^*_-$
pass to corresponding forms
$$
\langle\cdot,\cdot\rangle_+
\colon
g_+
\otimes
g_+
@>>>
\bold R,
\quad
\langle\cdot,\cdot\rangle_-
\colon
g_-
\otimes
g_-
@>>>
\bold R,
\quad
\langle\cdot,\cdot\rangle_0
\colon
g_0
\otimes
g_0
@>>>
\bold R;
\tag3.4
$$
and the orthogonal structure
(3.2) decomposes as
$$
\langle\cdot,\cdot\rangle'_{g}
=\langle\cdot,\cdot\rangle_+
+
\langle\cdot,\cdot\rangle_-
+
\langle\cdot,\cdot\rangle_0
\colon
g
\otimes
g
@>>>
\bold R.
\tag3.5
$$
Finally,
the 2-form
$\langle\cdot,\cdot\rangle_{g^*}$
on $g^*$ passes to
the adjoint action invariant
symmetric
bilinear form
$$
\langle\cdot,\cdot\rangle_{g}
=
\langle\cdot,\cdot\rangle_+
+
\left(-
\langle\cdot,\cdot\rangle_-\right)
$$
on $g$.
In other words,
all the relevant structure is preserved.
\smallskip
With reference to the decomposition (3.3),
the Lie algebra bundle
$\roman{ad}(\xi)$
decomposes
into  a direct sum
of corresponding Lie algebra bundles
$\zeta_+$,\ $\zeta_-$, and $\zeta_0$.
Hence
the graded vector space
$\Omega^*(\Sigma,\roman{ad}(\xi))$
of
$\roman{ad}(\xi)$-valued forms
decomposes as well into the direct sum of
$\Omega^*(\Sigma,\zeta_+)$,
$\Omega^*(\Sigma,\zeta_-)$,
and
$\Omega^*(\Sigma,\zeta_0)$
and, whatever connection $\widetilde A$, the operator
of covariant derivative
$d_{\widetilde A}$
preserves the decompositions.
When we divide out the appropriate groups of translations,
we obtain the affine spaces
$$
\gathered
\Cal A(\xi)_+ =\Cal A(\xi)\big /\left(
\Omega^*(\Sigma,\zeta_-)
\oplus
\Omega^*(\Sigma,\zeta_0)\right),
\\
\Cal A(\xi)_- =\Cal A(\xi)\big /\left(
\Omega^*(\Sigma,\zeta_+)
\oplus
\Omega^*(\Sigma,\zeta_0)\right),
\\
\Cal A(\xi)_0 =
\Cal A(\xi)\big /\left(
\Omega^*(\Sigma,\zeta_+)
\oplus
\Omega^*(\Sigma,\zeta_-)\right),
\endgathered
$$
together with a canonical isomorphism
of affine spaces
from
$\Cal A(\xi)$
to
\linebreak
$\Cal A(\xi)_+
\times
\Cal A(\xi)_-
\times
\Cal A(\xi)_0$.
The decomposition
(3.5)
of the orthogonal structure on $g$
entails that the relevant structure
made explicit in Section 1 of our paper
\cite \singula\
decomposes accordingly.
In particular,
the Lie algebra $\bold g(\xi)=
\Omega^0(\Sigma,\roman{ad}(\xi))$
of infinitesimal gauge transformations decomposes into a direct sum
of Lie ideals
$\bold g(\xi)_+=\Omega^0(\Sigma,\zeta_+)$,
\
$\bold g(\xi)_-=\Omega^0(\Sigma,\zeta_-)$,
and
$\bold g(\xi)_0=\Omega^0(\Sigma,\zeta_0)$.
\smallskip
We now consider the Yang-Mills theory
on $\xi$,
with the orthogonal structure
(3.2)
playing the role
of the orthogonal structure
on $g$ at the beginning of this Section.
The same
construction as
that of the weakly symplectic structure
$\sigma$
on $\Cal A(\xi)$
resulting from the latter
yields
weakly symplectic structures
$\sigma_+,\,
\sigma_-,\,
\sigma_0,\,
$
on
respectively
$\Cal A(\xi)_+,\,
\Cal A(\xi)_-,\,
\Cal A(\xi)_0,\,
$
resulting from the corresponding orthogonal structures (3.4).
Moreover the space
$\Cal N(\xi)\subseteq \Cal A(\xi)$ of central Yang-Mills connections
accordingly decomposes into a product
$
\Cal N(\xi)_+
\times
\Cal N(\xi)_-
\times
\Cal N(\xi)_0.
$
The results
in our paper
\cite\singula\
imply that,
near a point of
the top stratum
$N^{\roman{top}}(\xi)$,
the space
$N(\xi)$
locally decomposes into a product of
suitable submanifolds of
$ N(\xi)_+$,
$ N(\xi)_- $, and
$ N(\xi)_0$.
On taking on
these spaces
locally the Poisson structures
coming from
respectively
$\sigma_+$,
$-\sigma_- $,
and
the zero structure on
$ N(\xi)_0$,
we locally recover
the bracket
(2.1.3)
for an arbitrary 2-form $\langle\cdot,\cdot\rangle_{g^*}$.
This implies that the bracket (2.1.3) satisfies the
Jacobi identity in the general case.
\smallskip
When the form
$\langle\cdot,\cdot\rangle_{g^*}$
is non-degenerate,
in the above discussion, the objects labelled
$-_0$ do not occur. This implies that the resulting Poisson structure
(2.1.3) is then symplectic. The proof of Theorem 2.1 is now complete.
\smallskip\noindent
{\smc Remark 3.6.}
The argument may be extended to show that,
when the point $[A]$ of $N(\xi)$ represented by
$A$ is non-singular,
the Kuranishi map yields
Darboux coordinates
for the Poisson structure
on
$N(\xi)$ near $[A]$.

\medskip\noindent{\bf 4. Stratified symplectic structures}
\smallskip\noindent
We remind the reader that the
spaces
$N(\xi)$
and
$\roman{Rep}_{\xi}(\Gamma,G)$
are stratified
by connected components of orbit types;
see Section 2 of our paper
\cite\singulat\ for details.

\proclaim{Theorem 4.1}
Suppose the given  coadjoint action invariant
symmetric bilinear form
on $g^*$
is non-degenerate.
Then the
Poisson bracket
{\rm (2.1.3)}
restricts to a symplectic
Poisson bracket
on each stratum.
In other words,
this Poisson structure
endows
the spaces
$N(\xi)$
and
$\roman{Rep}_{\xi}(\Gamma,G)$
with a structure of a
stratified symplectic space.
\endproclaim

\smallskip
\noindent
{\smc Remark 4.2.}
We note that, in the statement of the theorem,
there is
{\it no need\/}
to assume the 2-form to be positive definite.

\demo{Proof of {\rm (4.1)}}
Clearly this is a local statement,
and we can argue in terms of the local model
of a neighborhood
of a point $[A]$ of $N(\xi)$
given in  \cite\singula\ (2.32),
with the following slight modification
which is needed when
the given 2-form
on $g^*$ is not positive definite:
Let $A$ be a central Yang-Mills connection,
and consider the  decomposition
of
$
\roman H^1_A(\Sigma,\roman{ad}(\xi))
$
into the direct sum of
$\roman H^1_A(\Sigma,\zeta_+)$
and
$\roman H^1_A(\Sigma,\zeta_-)$
where the notation is that
in the previous Section.
Denote by
$\sigma_A^+$ and $\sigma_A^-$
the symplectic structures
on
$\roman H^1_A(\Sigma,\zeta_+)$
and
$\roman H^1_A(\Sigma,\zeta_-)$
induced by
the 2-forms
$\langle\cdot,\cdot\rangle_+$
and
$\langle\cdot,\cdot\rangle_-$,
respectively,
cf. (3.4),
and
consider the mappings
$$
\Theta^+_A \colon
\roman H^1_A(\Sigma,\zeta_+)
@>>>
\roman H^2_A(\Sigma,\zeta_+),
\quad
\Theta^-_A \colon
\roman H^1_A(\Sigma,\zeta_-)
@>>>
\roman H^2_A(\Sigma,\zeta_-),
$$
given by the assignments
to
$\alpha \in \roman H^1_A(\Sigma,\zeta_+)$
and
$\beta \in \roman H^1_A(\Sigma,\zeta_-)$
of
$\Theta^+_A(\alpha) =\frac 12 [\alpha,\alpha]_A$ and
$\Theta^-_A(\beta) =- \frac 12 [\beta,\beta]_A$, respectively.
Their direct sum
$
\Theta_A
$
yields
a momentum mapping
for the
action of the stabilizer
$Z_A$
of $A$ in $\Cal G(\xi)$
on
$\roman H^1_A(\Sigma,\roman{ad}(\xi))$,
with the symplectic structure
$\sigma_A^+ + (-\sigma_A^-)$;
in fact it is the unique momentum mapping
having the value  zero at the origin.
Now the Marsden-Weinstein
reduced space
$
\roman H_A
=
\Theta_A^{-1}(0) \big / Z_A
$
is a local model of $N(\xi)$ near $[A]$,
cf. \cite\singula,\ \cite\smooth.
By the main result of
{\smc Sjamaar-Lerman}~\cite\sjamlerm,
the decomposition of
$\roman H_A$ into connected components
of orbit types
is a stratification in such a way that each stratum inherits
a symplectic structure.
Moreover, write
$\roman H^1_A
=
\roman H^1_A(\Sigma,\roman{ad}(\xi))$
for short, and
consider the algebra
$$
C^{\infty}(
\roman H_A)=
\left(C^{\infty}(\roman H^1_A)\right)^{Z_A}
\big / I_A^{Z_A}
$$
of smooth
$Z_A$-invariant functions
$\left(C^{\infty}(\roman H^1_A)\right)^{Z_A}$
on
$\roman H^1_A$
modulo its ideal
$I_A^{Z_A}$
of
functions
that vanish on
the zero locus
$\Theta_A^{-1}(0)$.
This algebra
endows
$\roman H_A$
with a smooth structure;
see \cite\smooth\ (6.1.2)
for more details.
By {\smc Arms-Cushman-Gotay}~\cite\armcusgo,
the symplectic Poisson
bracket on $C^{\infty}(\roman H^1_A)$
passes to a Poisson bracket $\{\cdot,\cdot\}_A$ on
$C^{\infty}(\roman H_A)$ which, in turn, restricts to
the symplectic
Poisson bracket on
each stratum of
$\roman H_A$.
This establishes the assertion of the Theorem locally.
\smallskip
By \cite\smooth\ (6.2),
near $[A] \in N(\xi)$,
with the smooth structure
$C^{\infty}(\roman H_A)$, the space
$\roman H_A$
is a local model of
$N(\xi)$, with its smooth structure
$C^{\infty}(N(\xi))$
introduced in Section 4 of \cite\smooth.
However, by construction,
near $[A] \in N(\xi)$,
the algebra $C^{\infty}(\roman H_A)$
with the
Poisson bracket
$\{\cdot,\cdot\}_A$
is in fact
a local model of
$N(\xi)$
with the Poisson bracket
on $C^{\infty}(N(\xi))$
induced from
(2.1.3)
via the Wilson loop mapping
$\rho_{\flat}$ from
$N(\xi)$ to $\roman{Rep}_{\xi}(\Gamma,G)$.
This completes the proof. \qed
\enddemo

The proof also establishes the following.

\proclaim{Addendum 4.3}
For every central Yang-Mills connection $A$,
near $[A] \in N(\xi)$,
the Poisson algebra $(C^{\infty}(\roman H_A),\{\cdot,\cdot\}_A)$
yields a local model of
$N(\xi)$
with its Poisson structure
and likewise,
near the point $\rho_{\flat}[A]$,
a local model of
$\roman{Rep}_{\xi}(\Gamma,G)$
with its Poisson structure,
where $\rho_{\flat}$
refers to the Wilson loop mapping
from
$N(\xi)$
to
$\roman{Rep}(\Gamma,G)$.
More precisely, the choice of $A$
(in its class $[A]$)
induces a Poisson diffeomorphism
of an open neighborhood
$W_A$ of $[0] \in \roman H_A$
onto an open neighborhood
$U_A$
of
$[A] \in N(\xi)$,
where
$W_A$ and
$U_A$
are endowed with the induced
smooth and Poisson structures,
and a similar statement
holds for
$\roman{Rep}_{\xi}(\Gamma,G)$
near the point $\rho_{\flat}[A]$.
\endproclaim

\medskip\noindent {\bf 5. Twist flows}\smallskip\noindent
Let $f$ be a smooth invariant real function
on $G$, and let $C$ be a closed curve in $\Sigma$,
having starting point $Q$. The homotopy class  of $C$ induces
a homomorphism
$[C]$ from
$\bold Z$ to the fundamental group $\pi$ of $\Sigma$, and the association
$\phi \mapsto f(\phi([C](1))) \in \bold R$
induces a real valued function $f^C$ on
$\roman{Rep}(\pi,G) = \roman{Hom}(\pi,G) \big / G$.
Let
$\xi$ be a flat $G$-bundle over $\Sigma$.
Since
$[C]$ lifts to a homomorphism
from
$\bold Z$ to the free group $F$
on the chosen generators,
$f^C$ yields
a function in $C^{\infty}(\roman{Rep}_\xi(\pi,G))$,
and hence $\{f^C,\cdot\}$ is a derivation of
$C^{\infty}(\roman{Rep}_\xi(\pi,G))$.
On each stratum it amounts of course to a smooth vector field.
The corresponding flow on the {\it non-singular\/} part
of $\roman{Rep}_\xi(\pi,G)$
has been studied by {\smc Goldman} in \cite\goldmtwo,
referred to as a {\it twist flow\/}.
However,
the derivation $\{f^C,\cdot\}$
in fact integrates to a
\lq\lq twist flow\rq\rq\  on the whole space
$\roman{Rep}_\xi(\pi,G)$,
that is, an action of the real line on this space
preserving the smooth structure.
The argument in Section 3 of \cite\sjamlerm,
applied to our local model constructed in \cite\singula,
shows that this is locally so, and using the partition of unity
established in our paper \cite\smooth,
we conclude that these twist flows in fact exist globally and are unique.
We hope to give the details at another occasion.

\bigskip
\centerline{\smc References}
\medskip
\widestnumber\key{999}

\ref \no  \armcusgo
\by J. M. Arms,  R. Cushman, and M. J. Gotay
\paper  A universal reduction procedure for Hamiltonian group actions
\paperinfo in: The geometry of Hamiltonian systems, T. Ratiu, ed.
\jour MSRI Publ.
\vol 20
\pages 33--51
\yr 1991
\publ Springer
\publaddr Berlin-Heidelberg-New York-Tokyo
\endref

\ref \no  \atibottw
\by M. Atiyah and R. Bott
\paper The Yang-Mills equations over Riemann surfaces
\jour Phil. Trans. R. Soc. London  A
\vol 308
\yr 1982
\pages  523--615
\endref

\ref \no  \goldmone
\by W. M. Goldman
\paper The symplectic nature of the fundamental group of surfaces
\jour Advances
\vol 54
\yr 1984
\pages 200--225
\endref

\ref \no  \goldmtwo
\by W. M. Goldman
\paper Invariant functions on Lie groups and Hamiltonian flows of
surface group representations
\jour Inventiones
\vol 85
\yr 1986
\pages 263--302
\endref

\ref \no \poisson
\by J. Huebschmann
\paper Poisson cohomology and quantization
\jour
J. f\"ur die Reine und Angewandte Mathematik
\vol  408
\yr 1990
\pages 57--113
\endref

\ref \no  \souriau
\by J. Huebschmann
\paper On the quantization of Poisson algebras
\book Symplectic Geometry and Mathematical Physics
\bookinfo Actes du colloque en l'honneur de Jean-Marie Souriau,
P. Donato, C. Duval, J. Elhadad, G.M. Tuynman, eds.;
Progress in Mathematics, Vol. 99
\publ Birkh\"auser
\publaddr Boston $\cdot$ Basel $\cdot$ Berlin
\yr 1991
\pages 204--233
\endref

\ref \no \singula
\by J. Huebschmann
\paper The singularities of Yang-Mills connections
for bundles
on a surface. I. The local model
\paperinfo Math. Z. (to appear), dg-ga/9411006
\endref
\ref \no  \singulat
\by J. Huebschmann
\paper The singularities of Yang-Mills connections
for bundles on a surface. II. The stratification
\paperinfo Math. Z. (to appear), dg-ga/9411007
\endref
\ref \no \topology
\by J. Huebschmann
\paper
Holonomies of Yang-Mills connections
for bundles  on a surface
with disconnected structure group
\jour Math. Proc. Cambr. Phil. Soc
\vol 116
\yr 1994
\pages 375--384
\endref

\ref \no \smooth
\by J. Huebschmann
\paper
Smooth structures on certain
moduli spaces for bundles on a surface
\paperinfo preprint 1992, dg-ga/9411008
\endref

\ref \no \singulth
\by J. Huebschmann
\paper The singularities of Yang-Mills connections
for bundles on a surface. III. The identification of the strata
\paperinfo in preparation
\endref

\ref \no \locpois
\by J. Huebschmann
\paper Poisson geometry of
flat connections for {\rm SU(2)}-bundles on surfaces
\paperinfo Math. Z. (to appear), hep-th/9312113
\endref

\ref \no \modusym
\by J. Huebschmann
\paper Symplectic and Poisson structures of certain moduli spaces
\paperinfo Preprint 1993, hep-th/9312112
\endref

\ref \no \huebjeff
\by J. Huebschmann and L. Jeffrey
\paper Group cohomology construction
of symplectic forms on certain moduli spaces
\jour Int. Math. Research Notices
\vol 6
\yr 1994
\pages 245--249
\endref

\ref \no \karshone
\by Y. Karshon
\paper
An algebraic proof for the symplectic
structure of moduli space
\jour Proc. Amer. Math. Soc.
\vol 116
\yr 1992
\pages 591--605
\endref

\ref \no \marswein
\by J. Marsden and A. Weinstein
\paper Reduction of symplectic manifolds with symmetries
\jour Rep. on Math. Phys.
\vol 5
\yr 1974
\pages 121--130
\endref

\ref \no \narashed
\by M. S. Narasimhan and C. S. Seshadri
\paper Stable and unitary vector bundles on a compact Riemann surface
\jour Ann. of Math.
\vol 82
\yr 1965
\pages  540--567
\endref

\ref \no \naramntw
\by M. S. Narasimhan and S. Ramanan
\paper Moduli of vector bundles on a compact Riemann surface
\jour Ann. of Math.
\vol 89
\yr 1969
\pages  19--51
\endref

\ref \no \naramnth
\by M. S. Narasimhan and S. Ramanan
\paper 2$\theta$-linear systems on abelian varieties
\jour Bombay colloquium
\yr 1985
\pages  415--427
\endref

\ref \no \seshaone
\by C. S. Seshadri
\paper Spaces of unitary vector bundles on a compact Riemann surface
\jour Ann. of Math.
\vol 85
\yr 1967
\pages 303--336
\endref

\ref \no \sjamlerm
\by R. Sjamaar and E. Lerman
\paper Stratified symplectic spaces and reduction
\jour Ann. of Math.
\vol 134
\yr 1991
\pages 375--422
\endref

\ref \no \weinstwe
\by A. Weinstein
\paper On the symplectic structure of moduli space
\paperinfo A. Floer memorial, Birkh\"auser Verlag, to appear
\endref

\ref \no  \whitnone
\by H. Whitney
\paper Analytic extensions of differentiable functions defined
on closed sets
\jour Trans. Amer. Math. Soc.
\vol 36
\yr 1934
\pages  63--89
\endref

\enddocument